\documentclass[UTF-8, 11pt, twoside, reqno]{amsart}
\usepackage[%colorlinks=false,   
colorlinks=true, linkcolor=blue, citecolor=magenta]{hyperref}

\usepackage{enumerate,enumitem,amsmath,amssymb,amsthm,amssymb,amsfonts,mathrsfs}
\usepackage{bbm}
\usepackage{cases}
\usepackage{graphicx}
\usepackage{stmaryrd}
\usepackage{mathptmx}
\usepackage{color}
\definecolor{MyDarkBlue}{cmyk}{0.8,0.3,0.8,0.4}
\definecolor{yellow}{rgb}{0.99,0.99,0.70}
\definecolor{white}{rgb}{1.0,1.0,1.0}
\definecolor{black}{rgb}{0.00,0.00,0.00}

\usepackage{fancyhdr}
{}
\makeatletter
\newcommand{\mylabel}[2]{#2\def\@currentlabel{#2}\label{#1}}
\makeatother

\numberwithin{equation}{section}

\newcommand{\be}{\begin{eqnarray}}
\newcommand{\ee}{\end{eqnarray}}
\newcommand{\ce}{\begin{eqnarray*}}
\newcommand{\de}{\end{eqnarray*}}
\newtheorem{theorem}{Theorem}[section]
\newtheorem{lemma}[theorem]{Lemma}
\newtheorem{remark}[theorem]{Remark}
\newtheorem{definition}[theorem]{Definition}
\newtheorem{proposition}[theorem]{Proposition}
\newtheorem{Examples}[theorem]{Example}
\newtheorem{corollary}[theorem]{Corollary}

\def\eps{\varepsilon}

\def\e{\mathrm{e}}

\def\a{\alpha}
\def\om{\omega}
\def\Om{\Omega}

\def\p{\partial}
\def\d{\delta}

\def\[{{\Big[}}
\def\]{{\Big]}}
\def\<{{\langle}}
\def\>{{\rangle}}
\def\({{\Big(}}
\def\){{\Big)}}

\def\bx{{\mathbf{x}}}

\def\dif{{\mathord{{\rm d}}}}

\def\no{\nonumber}
\def\={&\!\!=\!\!&}

\def\cJ{{\mathcal J}}

\def\cL{{\mathcal L}}

\def\cN{{\mathcal N}}

\def\sS(\mR^d){{\mathcal S}}
\def\cT{{\mathcal T}}

\def\mD{{\mathbb D}}

\def\mI{{\mathbb I}}

\def\mN{{\mathbb N}}

\def\mR{{\mathbb R}}

\def\mZ{{\mathbb Z}}

\def\bP{{\mathbf P}}

\def\bE{{\mathbf E}}
\def\1{{\mathbf{1}}}

\def\sF{{\mathscr F}}

\def\sS{{\mathscr S}}

\def\sW{{\mathscr W}}

\def\E{\mathbf E}

\def\geq{\geqslant}
\def\leq{\leqslant}

\def\div{\mathord{{\rm div}}}
\def\iint{\int\!\!\!\int}

\def\eps{\varepsilon}

\def\e{\mathrm{e}}

\def\a{\alpha}
\def\om{\omega}
\def\Om{\Omega}

\def\p{\partial}
\def\d{\delta}

\def\[{{\Big[}}
\def\]{{\Big]}}
\def\<{{\langle}}
\def\>{{\rangle}}
\def\({{\Big(}}
\def\){{\Big)}}

\def\bx{{\mathbf{x}}}

\def\dif{{\mathord{{\rm d}}}}

\def\no{\nonumber}
\def\={&\!\!=\!\!&}
\def\bt{\begin{theorem}}
\def\et{\end{theorem}}
\def\bl{\begin{lemma}}
\def\el{\end{lemma}}
\def\br{\begin{remark}}
\def\er{\end{remark}}
\def\bx{\begin{Examples}}
\def\ex{\end{Examples}}
\def\bd{\begin{definition}}
\def\ed{\end{definition}}
\def\bp{\begin{proposition}}
\def\ep{\end{proposition}}
\def\bc{\begin{corollary}}
\def\ec{\end{corollary}}

\def\geq{\geqslant}
\def\leq{\leqslant}

\def\div{\mathord{{\rm div}}}
\def\iint{\int\!\!\!\int}
\def\si{\sigma}
 \def\mR{\mathbb R}
 \def\mR{\mathbb R}    
\def\N{\mathbb N}  
   
\def\<{\langle} \def\>{\rangle}  
    
\def\d{\text{\rm{d}}}   
  \def\si{\sigma} 
 \def\beq{\begin{equation}}  
 
\def\e{\text{\rm{e}}}    
  
 \def\P{\mathbf P}

\allowdisplaybreaks

\begin{document}
%\pagecolor{MyDarkBlue}\color{white}

\title{Regularity properties of jump diffusions with irregular coefficients.}

\author{Guohuan Zhao} 

\address{Guohuan Zhao: Department of Mathematics, Bielefeld University, Germany \\
Email: zhaoguohuan@gmail.com
 }

\thanks{Research of Guohuan is supported by the German Research Foundation (DFG) through the Collaborative Research Centre(CRC) 1283 Taming uncertainty and profiting from randomness and low regularity in analysis, stochastics and their applications.}

\begin{abstract}
In this paper, we investigate the regularity properties of strong solutions to SDEs driven by L\'evy processes with irregular drift coefficients. Under some mild conditions, we show that the singular SDE has a unique strong solution for each starting point and the family of all the solutions forms a stochastic flow. Moreover, the Malliavin differentiability of the strong solutions is also obtained. As an application, we also prove a path-by-path uniqueness result for some related random ODEs. Our proofs are based on Schauder's estimate for the corresponding resolvent equation, which is also analytically interesting in itself. 
\end{abstract}

\maketitle



\section{Introduction}
\subsection{Motivation and main results}
The motivation of this paper is to study the following stochastic differential equation (SDE) driven by L\'evy process with irregular drift coefficient  
\begin{equation}\label{Eq-SDE}
\d X_t(x)=b(X_t(x))\d t+\sigma (X_{t-}(x))\d Z_t, ~~X_{0}(x)=x\in\mR^d.  
\end{equation}
Here $Z$ is a $d$-dimensional pure jump L\'evy process 
with L\'evy-Khintchine triplet $(0,0,\nu)$, 
$b:\mR^d\to\mR^d$ and $\sigma:\mR^d\to\mR^d\otimes\mR^d$ are measurable maps. Two basic questions about SDEs with jumps as in \eqref{Eq-SDE} are its (strong) well-posedness and the regularity properties of the solutions $X_t(x, \om)$ with respect to (w.r.t.) the initial value $x$ and the sample point $\om$. When the coefficients are regular, the reader is referred to Kunita's fundamental work \cite{kunita2004stochastic} for the existence and regularity of stochastic flows associated to SDEs driven by general Poisson random measures. This paper mainly focuses on the case when the drift coefficients are irregular.  
In the case that $Z_t$ is a truncated rotational symmetric $\a$-stable process with $\a \in (1, 2)$, $\si=\mI$ and $b\in C^\beta_b$ with $\beta>2-\a$, Haadem and Proske \cite{haadem2014construction} proved the existence and Malliavin differentiability of $X_{t}(x,\om)$. Another aim of this paper is to extend Haadem and Proske's results to more general (truncated) $\a$-stable like process $Z$ and relax the restrictions on $\a$, $\beta$ at the same time. 

To achieve our goals, in this paper we first use a Zvonkin's change of variables to remove the irregular drift terms and transform the old SDEs to some new SDEs with regular coefficients, then we prove the desired properties of the solutions to \eqref{Eq-SDE} by studying the corresponding properties of the solutions to the new equation. From the probabilistic point of view, this  approach can not be considered novel and the main idea can be traced back at least to \cite{zvonkin1974transformation}. However, in order to get the so-called Zvonkin's transform, we need to resort to some nice apriori estimates for the following resolvent equation associated to \eqref{Eq-SDE} 
\begin{align}\label{Eq-PDE} 
\lambda u-\cL^b u=f, 
\end{align}
where $\cL^b$ is a L\'evy type operator defined by  
\begin{align*}
\cL^bu(x):=& \cL u(x)+ b(x)\cdot\nabla u(x)\\
:=& \int_{\mR^d} \left[u(x+\sigma(x)z)-u(x)-\nabla u(x)\cdot \sigma(x) z\1_{B_1}(z)\right]\nu(\dif z) + b(x)\cdot\nabla u(x), 
\end{align*}
and $\nu$ is the L\'evy measure of the $\a$-stable like process $Z$. Noting that when $\a\in (0,1)$, the gradient $\nabla$ is of higher order than the nonlocal dissipative operator $\cL$, the regularity theory for \eqref{Eq-PDE} is a highly nontrivial analytic problem when $\a\in (0,1)$. 

To state our results precisely,  we start by introducing the assumptions needed in this paper. Let $\a\in (0,2)$ and $Z$ be a L\'evy process with L\'evy-Khintchine triplet $(0,0,\nu)$.  
\begin{enumerate}[leftmargin=1.05cm]
\item [\mylabel{A_nu}{{\rm (A$_1$)}}] 
$\nu$ is symmetric and there are some constants $c_0, \rho\in (0,1)$ such that 
\be\label{Eq-nu-lower} 
\int_{B_r} \<\theta, z\>^2 \nu(\d z)\geq c_0 r^{2-\a}, \quad \forall r\in (0,\rho), \  \theta\in \{\vartheta\in \mR^d: |\vartheta|=1\},
\ee
and 
\be\label{nu-upper}
\int_{B_r} |z|^2 \nu(\d z)\leq c_0^{-1}r^{2-\a}, \quad \forall r\in (0,\rho), 
\ee
\item [\mylabel{A_si}{{\rm (A$_2$)}}] $\si$ is non-degenerate and there is a constant $\Lambda\geq 1$ such that
$$
\|\nabla \sigma\|_\infty \leq \Lambda \text{\quad and \quad} \Lambda^{-1} |\xi|\leq |\sigma(x) \xi|\leq \Lambda |\xi|,
$$
\item [\mylabel{A_b}{{\rm (A$_3$)}}] $b\in C^\beta_b$ with $\beta>1-\frac{\a}{2}$, 
\item [\mylabel{A_nu2}{{\rm (A$_{4}$)}}]  
$\si\in C_b^{1+\delta}$ for some $\delta>0$. $\nu$ has a compact support and  
$$
\mathrm{supp}\, \nu\subseteq B_{r_0} \text{ with } r_0=\|\nabla \si\|_\infty^{-1}. 
$$
\end{enumerate}

The main analytic result of this paper reads as follows.  
\bt\label{Th-Holder}
Let $\alpha\in(0,2)$ and $(1-\a)^+<\gamma<\beta<1$.  Assume $\nu$ is supported on $B_R$ satisfying conditions \ref{A_nu}, $\si$ satisfies \ref{A_si} and $b\in C_b^\beta$. Then there is a constant $\lambda_0$ such that for any $\lambda\geq \lambda_0$ and $f\in C_b^\beta$, equation \eqref{Eq-Resolvent} has a unique solution in $C_b^{\a+\gamma}$. Moreover, 
\begin{align}\label{Eq-Holder-Est}
\lambda\|u\|_{C_b^\gamma}+\|u\|_{C_b^{\alpha+\gamma}} \leq C \|f\|_{C_b^\beta}, 
\end{align}
where the constants $\lambda_0, C$ only depend on $d, \alpha,\beta,\gamma, R, c_0, \rho, \Lambda $ and $\|b\|_{C_b^{\beta}}$.   
\et

\br\label{Rek-Const}
\begin{enumerate}
\item The assumption on $\si$ in Theorem \ref{Th-Holder} can be weakened, but we do not attempt to do that here, since the above results are enough for our main purpose. 
\item Suppose $\si$ is independent of $x$. Assume $\nu$ only satisfies \ref{A_nu} (does not need to be compactly supported) and $b\in C_b^\beta$ with $\beta\in (1-\frac{\a}{2},1)$. Then by Lemma \ref{Priori1} below and the last part of proof for Theorem \ref{Th-Lp-Est}, one can see that for any $\gamma\in ((1-\a)^+, \beta)$ there is a constant $\lambda_0$ such that for any $\lambda\geq \lambda_0$ and $f\in C_b^\beta$, equation \eqref{Eq-Resolvent} has a solution $u\in C_b^{\a+\gamma}$. Moreover, the estimate  \eqref{Eq-Holder-Est} still holds. 
\end{enumerate}
\er

Thanks to the mentioned Zvonkin's transformation, the above result can be applied to prove the regularity properties of solutions to stochastic differential equation \eqref{Eq-SDE}. 
 
\bt\label{Th-Main}
Under conditions \ref{A_nu}, \ref{A_si} and \ref{A_b}, for each starting point $x\in \mR^d$, there is a unique strong solution $X_t(x)$ to equation \eqref{Eq-SDE}. Moreover, 
\begin{enumerate}[leftmargin=0.9cm, label = (\roman*), ref=(\roman*)]
\item  [\mylabel{prop: i} {(i)}] if the jumping size of $Z_t$ is bounded, then for each $t\geq 0$, $X_t(x)$ is Malliavin differentiable; 
\item [\mylabel{prop: ii} {(ii)}]  if $\si$ and $\nu$ satisfy condition \ref{A_nu2} or $\si$ does not depend on $x$, then $\{X_t(x)\}_{t\geq0; x\in \mR^d}$ forms a $C^1$-stochastic diffeomorphism flow.
\end{enumerate}
\et 
\br
\begin{enumerate}[leftmargin=0.9cm]
\item The strong well-posedness of \eqref{Eq-SDE} has been essentially established in \cite{chen2021supercritical}, so the main contribution  of this work lies in Theorem \ref{Th-Main} \ref{prop: i} and \ref{prop: ii}. 
\item All the above results are 
uniform for all truncated non-degenerate symmetric $\a$-stable like processes with $\a\in (0,2)$ whose L\'evy measure can be singular w.r.t.  the Lebesgue measure on $\mR^d$ (one typical example of $Z$ is the truncated cylindrical $\a$-stable process). Moreover, our condition on $b$ is strictly weaker than that in \cite{haadem2014construction}. 
\item  For the existence of $C^1$-stochastic diffeomorphism flow, the compactly supported assumption on $\nu$ is not needed when $\sigma$ is independent of $x$. 
\end{enumerate}
\er

To summarize, this article offers the following novel contributions 
\begin{itemize}
\item in aspect of nonlocal partial differential equations, both $L^p$ type and Schauder type estimates are obtained for equation \eqref{Eq-PDE}, where the L\'evy measure $\nu$ associated to the L\'evy type operator $\cL$ can be singular w.r.t. the Lebesgue measure on $\mR^d$; 
\item in aspect of stochastic analysis, regularity of stochastic flows for jump SDEs with irregular coefficients driven by multiplicative noises is studied under very mild assumptions. 
\end{itemize}

\subsection{Background and main approach}
In the past two decades, there is a great interest for both analysts and probabilists to study nonlocal L\'evy type operators and corresponding Markov processes. Parts of the reasons lie in the fact that they have many applications in mathematical finance, control, physics, image processing, etc., and have various connections with other branches of mathematics. Interested readers are referred to  \cite{bottcher2013levy, jacob2001levy, jacob2005pseudo, schilling2009symbol} and the references therein for the connection between  pseudo-differential operators and L\'evy type processes,  \cite{bass2002transition, bass2002harnack,  bass2005holder, chen2003heat}, etc. for the study of potential theory of nonlocal operators via probabilistic approach and \cite{caffarelli2011regularity, kassmann2009priori, silvestre2006holder}, etc. for H\"older regularity estimates of solutions to nonlocal PDEs.  

For the regularity estimates for nonlocal PDEs with first order terms, to our best knowledge, it was Silvestre \cite{silvestre2012differentiability} who first  obtained a Schauder type estimate for the parabolic analogy of \eqref{Eq-PDE} when $\cL$ is the usual fractional Laplacian $\Delta^{\alpha/2}:=-(-\Delta)^{\alpha/2}$ with $\alpha\in(0,1)$ and $b\in C_b^\beta$ with $\beta\in(1-\alpha,1)$. Around the same time, \eqref{Eq-PDE} was also studied by Priola  \cite{priola2012pathwise, priola2015stochastic} in H\"older spaces for $\a$-stable like operators with $\a\in [1,2)$. The reader can also refer to \cite{de2020schauder, dong2018dini, dong2013schauder, hao2020schauder, lin2019nonlocal} for some important extensions and \cite{zhang2018dirichlet} for the Dirichlet problem associated with $\cL+b\cdot\nabla $ in the parabolic setting. However, we emphasize that all the above literatures can not cover the important case $\cL=\sum_{i=1}^d -(-\p_{ii}^2)^{\a/2}$
with $\a\in (0,1/2]$. In this paper, we give a systematic and self-contained approach to the apriori estimates for the solutions of resolvent equations, which allows one to handle a large class of L\'evy type operator in a uniform way, in particular, for operators associated with cylindrical $\a$-stable processes with $\a\in (0,2)$. Since $\cL$ is not the dominant term when $\a\in (0,1)$, it is very hard to apply the perturbation arguments as in the case $\a\in (1,2)$. To overcome this difficulty, we resort to some techniques in harmonic analysis, for instance, Littlewood-Paley decomposition and Bernstein type inequalities. To be specific, when $\si$ is a matrix, we first give a key Bernstein type  inequality \eqref{Eq-Bernstein2} with a self-contained proof.  Then we use the energy method to establish an apriori estimate for the solution of \eqref{Eq-PDE} in the Sobolev-Slobodeckij  space $W^s_p$, see Lemma \ref{Priori1}. For the  general case with varying coefficient $\sigma$,  under the additional condition that $\nu$ is compactly supported, we also obtain the corresponding conclusion by the usual freezing coefficients method. Using a localization technique from \cite{zhang2020stochastic} and \cite{xia2020lqlp},  we then get the desired Schauder type estimate. We need to point out that in  \cite{chen2021supercritical}, Chen, Zhang and the author of this paper adopted the similar approach of this paper to study the Kolmogorov equation in Besov space $B^s_{p,\infty}$. But the main analytic result (Theorem 3.6) therein needs a special assumption on the oscillation of $\si$, since $B^s_{p,\infty}$ does not enjoy the localization principle (see Lemma \ref{Le-Partition} below). This prevents us from applying their results to our  probabilistic purpose.

\medskip

One the other hand, the study of SDEs  driven by non-degenerate Gaussian noises with irregular coefficients has a long history and the literature on this topic is quite impressive. After Zvonkin's pioneering work \cite{zvonkin1974transformation}, there have  been a lot of works devoted to studying the strong well-posedness of singular SDEs, among which, we quote \cite{krylov2005strong, veretennikov1980strong2, zhang2011stochastic, zhang2016stochastic}. Moreover, there are also many articles investigating the regularity properties of $X_t(x,\om)$ w.r.t.  $x$ and $\om$. Among all, let us mention that the weak differentiability of $X_t(x)$ w.r.t. $x$ in the multiplicative noise case was shown in \cite{zhang2011stochastic, zhang2016stochastic} by using  Zvonkin's idea. In \cite{menoukeu2013variational, meyer2010construction, mohammed2015sobolev}, the authors adopted another approach based on Malliavin calculus to study the solvability of \eqref{Eq-SDE} as well as the Malliavin differentiability of the solutions in the additive noise case. 

When the driven noise $Z$  is an $\a$-stable process,  Tanaka, Tsuchiya and Watanabe \cite{tanaka1974perturbation} proved that if $d=1, \si=1$, $Z$ is a symmetric  $\alpha$-stable process with $\alpha\in[1,2)$ and $b$ is bounded measurable, then pathwise uniqueness holds for SDE \eqref{Eq-SDE}. They further show that when $\alpha\in(0,1)$,  even if $b$ is $\beta-$H\"older continuous with $\beta\in (0,1-\a)$, the pathwise uniqueness may fail. For multidimensional case, when $Z$ is an $\a$-stable like process with $\a\in [1,2)$ and $b$ is $\beta$-H\"older continuous  with $\beta>1- \alpha/{2}$, the well-posedness of \eqref{Eq-SDE} was studied by Priola \cite{priola2012pathwise, priola2015stochastic}, see also \cite{zhang2013stochastic} for the case when $b$ is in some Sobolev spaces. Later, the well-posedness in the supercritical case ($\a\in (0,1)$) was  studied in \cite{chen2018stochastic, chen2021supercritical}.  Weak differentiability of $X_t(x)$ w.r.t. $x$ together with the stochastic flow property were also discussed in \cite{chen2018stochastic}, when $\si=\mI$ and the semigroup of $Z$ fulfills some special gradient estimates. Using Malliavin's calculus for jump processes, Haadem and Proske \cite{haadem2014construction} studied the strong existence and  Malliavin differentiability by a similar approach as in \cite{menoukeu2013variational, meyer2010construction}. However, they had to assume that $Z_t$ is a truncated rotational symmetric $\a$-stable process with $\a>1$, $\si=\mI$ and $b\in C_b^\beta$ with $\beta>2-\a$, which are much stronger than our assumptions. For the weak well-posedness of \eqref{Eq-SDE}, the existence and continuity of the transition probability
density of the corresponding Markov process, the readers are referred to \cite{kulik2019weak, zhao2019weak}. We also note that \cite{huang2018euler, kuhn2019strong, mikulevivcius2018rate} study the convergence of the Euler-Maruyama approximation for \eqref{Eq-SDE}.

\medskip
 
The rest of this paper is organized as follows: In Section 2, we introduce some concepts and facts from Littlewood-Paley theory as well as Malliavin calculus for L\'evy processes. In Section 3, we study the nonlocal resolvent equation \eqref{Eq-PDE} and obtain some apriori estimates in both localized Sobolev spaces and H\"older spaces. In Section 4, we prove our probabilistic results by using our analytic results proved in Section 3 and Zvonkin's transform.

\section{Preliminary}
In this section, we first introduce the (localized) fractional Sobolev spaces and Besov spaces, which are the main function spaces to solve the PDE \eqref{Eq-PDE} in this paper.  For the sake of completeness,  we then present some concepts and results about the Malliavin calculus for L\'evy processes. 
\subsection{(Localized) Sobolev spaces and Besov spaces} Let $\chi:\mR^d\to[0,1]$ be a smooth radial function satisfying 
$$
\chi(\xi)=1,\ |\xi|\leq 1,\ \chi(\xi)=0,\ |\xi|\geq 3/2. 
$$
Define 
\be\label{Def-chiz}
\chi_z(x):= \chi(x-z).
\ee
\bd
For $n\in \mN$, $s\geq 0$ and $p\in (1,\infty)$, 
\begin{enumerate} 
\item the usual Sobolev space is defined as 
$$
W^n_p:= \left\{ f\in L^p: \|f\|_{W^n_p}:=\sum_{k=0}^n \|\nabla^k f\|_{p}<\infty\right\};  
$$
\item the Sobolev-Slobodeckij semi-norm is defined by
$$
{\displaystyle [f]_{\theta ,p}:=\left(\int\!\!\!\int _{\mR^d\times \mR^d }{\frac {|f(x)-f(y)|^{p}}{|x-y|^{\theta p+d}}}\;\d x\;dy\right)^{\frac {1}{p}}} ,  \  \theta \in (0,1) 
$$
and the Sobolev-Slobodeckij space $W^s_p$ is defined as 
$$
{\displaystyle W^{s}_{p}:=\left\{f\in W^{\lfloor s\rfloor }_p: \|f\|_{W^{s}_{p}}:=\|f\|_{W^{\lfloor s\rfloor}_p}+\sup _{|\alpha |=\lfloor s\rfloor }[\p^{\alpha }f]_{s-\lfloor s\rfloor ,p }<\infty \right\}}, \ s\notin \mN;  
$$
\item the usual Bessel potential space $H^s_p$ is defined as 
$$
H^s_p:=\left\{ f\in L^p: \|f\|_{H^s_p}:=\|(\mI-\Delta)^{s/2} f\|_p<\infty\right\}.   %\asymp \|f\|_p+\|(-\Delta)^{s/2}f\|_p. 
$$
\end{enumerate}
\ed
Now we can define the localized fractional Sobolev space, which will play a crucial role in our study of \eqref{Eq-PDE}. 
\bd 
For  $s\geq 0$, $p\in (1,\infty)$, define  
$$
\sW^s_p:=\left\{f\in W^s_{p,loc}: \|f\|_{\sW^s_p}:=\sup_{z\in \mR^d} \|f\chi_z\|_{W^s_p}<\infty\right\}.  
$$
\ed
%One can easily verify that $\sW^s_p$ is a Banach space. 
Next we recall some basic facts of the  Littlewood-Paley theory and give the definition of Besov spaces. Let $\sS(\mR^d)$ be the Schwartz space of all rapidly decreasing functions, and $\sS'(\mR^d)$ be the dual space of $\sS(\mR^d)$ 
called Schwartz generalized function (or tempered distribution) space. Given $f\in\sS(\mR^d)$,
let $\sF f=\hat f$  be the Fourier transform of $f$ defined by
$$
\hat f(\xi):=(2\pi)^{-d/2}\int_{\mR^d}\e^{-\mathrm{i}\xi\cdot x} f(x)\dif x.
$$
Define
$$
\varphi(\xi):=\chi(\xi)-\chi(2\xi).
$$
It is easy to see that $\varphi\geq 0$ and supp $\varphi\subset B_{3/2}\setminus B_{1/2}$ and
\begin{align*}
\chi(2\xi)+\sum_{j=0}^k\varphi(2^{-j}\xi)=\chi(2^{-k}\xi)\stackrel{k\to\infty}{\to} 1.
\end{align*}
In particular, if $|j-j'|\geq 2$, then
$$
\mathrm{supp}\varphi(2^{-j}\cdot)\cap\mathrm{supp}\varphi(2^{-j'}\cdot)=\varnothing.
$$
Let $h:=\sF^{-1} \chi$ be the inverse Fourier transform of $\chi$. Define
$$
h_{-1}(x):=\sF^{-1} \chi(2\cdot)(x)=2^{-d}h(2^{-1}x)\in\sS(\mR^d),
$$
and for $j\geq 0$,
\begin{align}\label{Eq-hj}
h_j(x):=\sF^{-1}\varphi(2^{-j}\cdot)(x)=2^{jd}h(2^jx)-2^{(j-1)d}h(2^{j-1}x)\in \sS(\mR^d).
\end{align}
Operators $\Delta_j$ and $S_j$ is defined by 
\begin{align*}
\Delta_j f:=
\left\{
\begin{aligned}
\sF^{-1}(\chi(2\cdot) \sF f), & j=-1, \\
\sF^{-1}(\varphi(2^{-j}\cdot) \sF f),& j\geq 0,
\end{aligned}
\right. \quad S_j f:=  \sum_{k\leq j-1} \Delta_k f. 
\end{align*}

It is easy to see that 
\begin{align}\label{Eq-conv-delta}
\Delta_j f(x)=(h_j*f)(x)=\int_{\mR^d}h_j(x-y)f(y)\dif y,\ \ j\geq -1.
\end{align}
\bd
For $s\in\mR$ and $p,q\in[1,\infty]$, the Besov space $B^s_{p,q}$ is defined as the set of all $f\in\sS'(\mR^d)$ with
$$
\|f\|_{B^s_{p,q}}:=\left(\sum_{j\geq -1}2^{jsq}\|\Delta_j f\|_p^q\right)^{1/q}<\infty. 
$$
Assume that $k\in \mN$, $\beta\in [0,1)$ and $s=k+\a$, %{\red with a slight abuse of notation,} 
let $C^s_b$ denote the collection of all bounded functions such that 
$$
\|f\|_{C_b^s}:= \sum_{i=0}^{k}\|\nabla^k f\|_\infty+ \sup_{x\neq y} \frac{|\nabla^k f(x)-\nabla^k f(y)|}{|x-y|^\beta}<\infty. 
$$
\ed

Last, we present two Bernstein type inequalities and the interrelationships of the above function  spaces. The following two lemmas can be found in \cite{triebel1992theory}. 

\bl\label{Le-BI}

Let $1\leq p\leq q\leq \infty$ and $j\geq 0$.  For any $f\in \sS'(\mR^d)$ with $\Delta_j f\in L^p$, we have
\begin{align}\label{Eq-BI1}
\|\nabla^k\Delta_j f\|_q\leq C_p 2^{(k+d(\frac{1}{p}-\frac{1}{q}))j}\|\Delta_jf\|_p,\  k=0,1,\cdots,
\end{align}
and
\begin{align}\label{Eq-BI2}
\|(-\Delta)^{s/2}\Delta_j f\|_q\leq C_p 2^{(s+d(\frac{1}{p}-\frac{1}{q}))j}\|\Delta_jf\|_p,\   s\in  \mR.  
\end{align}
\el
\bl
Let $s>\eps>0$ and $p\in (1,\infty)$.   
\begin{enumerate}[leftmargin=0.85cm]
\item 
It holds that $H^{s+\eps}_p\subseteq W^s_p\subseteq H^{s-\eps}_p$. 
\item If $0<s-\frac{d}{p}\notin \mN$, then $H^{s}_p\subseteq C_b^{s-\frac{d}{p}}$ and $W^{s}_p\subseteq C_b^{s-\frac{d}{p}}$.  
\item If $0<s\notin \N$, then $B^s_{\infty,\infty}= C_b^s$ and  $B^{s}_{p,p} = W^s_p$. 
\end{enumerate}
\el
The following lemma reveals the relationship between $\sW^s_p$ and $C_b^s$. 
\bl\label{Le-embed}
\begin{enumerate}[leftmargin=0.85cm]
\item For any  $\gamma>d/p$, we have $\sW^\gamma_p\subseteq C_b^{\gamma-\frac{d}{p}}$; 
\item For any $\beta>\gamma\geq 0$, we have  $C_b^\beta\subseteq \sW^\gamma_p$. 
\end{enumerate}
\el
\begin{proof}
The first conclusion is just a consequence of Sobolev embedding theorem. We only need to prove the second conclusion when $0\leq \gamma<\beta<1$. Obviously, if $u\in C_b^\beta$, then 
$$
\|u\chi_z\|_p\lesssim \|u\|_{L^\infty}\lesssim \|u\|_{C_b^\beta}, 
$$
where $\chi_z$ is defined by \eqref{Def-chiz}. 
Hence, $C_b^\beta\subseteq \sW^0_p$.  Note  $u\chi_z(x)=u\chi_z(y)=0$ if $y\in B_{1}(x)$ and $x\notin B_{5/2}(z)$. So for any $0< \gamma<\beta<1$, we have 
\begin{align*}
&\int\!\!\!\int_{\mR^d\times\mR^d}\frac{|u\chi_z (x)-u\chi_z(y)|^p}{|x-y|^{d+\gamma p}} \d x\d y\\
\lesssim & \int\!\!\!\int_{|x-z|\leq \frac{5}{2},\, |x-y|\leq 1} \frac{|u\chi_z (x)-u\chi_z(y)|^p}{|x-y|^{d+\gamma p}} \d x\d y+ \int\!\!\!\int_{|x-y|>1} \frac{|u\chi_z (x)-u\chi_z(y)|^p}{|x-y|^{d+\gamma p}} \d x\d y\\
\lesssim & \int\!\!\!\int_{|x-z|\leq \frac{5}{2},\, |x-y|\leq 1} \frac{|\chi_z (x)|^p \cdot |u(x)-u(y)|^p}{|x-y|^{d+\gamma p}} \d x\d y\\
&+\int\!\!\!\int_{|x-z|\leq \frac{5}{2},\, |x-y|\leq 1} \frac{|u (y)|^p \cdot |\chi_z(x)-\chi_z(y)|^p}{|x-y|^{d+\gamma p}} \d x\d y \\
&+  \int\!\!\!\int_{|x-y|>1} \frac{|u\chi_z (x)|^p+|u\chi_z(y)|^p}{|x-y|^{d+\gamma p}} \d x\d y\\
\lesssim &  \int\!\!\!\int_{|x-z|\leq \frac{5}{2},\, |x-y|\leq 1} \frac{\|u\|_{C_b^\beta}^p |x-y|^{\beta p}}{|x-y|^{d+\gamma p}} \d x\d y+ \int\!\!\!\int_{|x-y|> 1} \frac{|u\chi_z(x)|^p}{|x-y|^{d+\gamma p}} \d x\d y\\ 
\lesssim & \|u\|_{C_b^\beta}^p \int_{|x-z|\leq \frac{5}{2}} \d x\int_{|w|\leq 1} |w|^{-d+(\beta-\gamma)p} \d w+ \|u\|_{L^\infty}^p \int_{|x-z|\leq \frac{3}{2}} \d x\int_{|w|>1} |w|^{-d-\gamma p} \d w\\
\lesssim& \|u\|_{C_b^\beta}^p. 
\end{align*}
This yields, $\sup_{z\in \mR^d}\|u\chi_z\|_{W^\gamma_p} \lesssim C \|u\|_{C_b^\beta}$.  So we complete our proof. 
\end{proof}

\subsection{Malliavin Derivate for L\'evy processes}
In this subsection, we introduce some basic concepts of Malliavin calculus for L\'evy processes. The reader can refer to \cite{di2009malliavin} for more details. 

Suppose $N(\d t, \d x)$ is a Poisson point process with intensity measure $\nu(\d z)$. Let $\{\mathcal{F}_t\}_{0\leq t\leq T}$ be the filtration generated by $N$ and $\widetilde{N}(\d t,\d z):=N(\d t,\d z)-\nu(\d z)\d t$. 

\medskip

For each $n\in \mN_+$ and $f\in L^2(([0,T]\times \mR^d)^n; (\lambda\times \nu)^n)$, define 
\be\label{Eq-def-tildef}
\widetilde f(t_1, z_1; \cdots; t_n, z_n):= \frac{1}{n!}  \sum_{\si\in S_n} f(t_{\si(1)}, z_{\si(1)}; \cdots; t_{\si(n)}, z_{\si(n)}). 
\ee
The collection of all square
integrable symmetric functions is denoted by $\widetilde L^2(([0,T]\times\mR^d)^n, (\lambda\times \nu)^n)$(abbreviated by $\widetilde L^2((\lambda\times \nu)^n)$). For any $n\in \mN_+$ and $f_n\in \widetilde L^2((\lambda\times \nu)^n)$, we set 
$$
I_n(f_n):= \int_{([0,T]\times \mR^d)^n} f_n(t_1,z_1;\cdots;t_n,z_n)\widetilde{N}^{\otimes n}(\d {\bf t},\d {\bf z}),  \ {\bf t}=(t_1,\cdots t_n), {\bf z}=(z_1,\cdots, z_n)
$$
and $I_0(f_0):=f_0$ for the constant values $f_0\in \mR$. 
\begin{definition}
The stochastic Sobolev space $\mathbb{D}^1_2$ consists of all  $\mathcal{F}_T$ measurable random variables $F\in L^2(\P)$ with chaos expansion $F=\sum_{n=0}^\infty I_n(f_n)$ satisfying 
$$
\sum_{n=1}^\infty n n!\|f_n\|_{L^2((\lambda\times\nu)^n)}^2<\infty. 
$$
Define the Malliavin derivative operator $D: \mD^1_2\ni F \mapsto DF \in L^2 (\lambda\times \nu\times \bP)$  
by 
$$
D_{t,z}F:= \sum_{n=1}^\infty n I_{n-1}(\widetilde f_n(\cdot; t,z)),  \quad F\in \mD^1_2, 
$$
where $\widetilde f_n(\cdot; t, z)$ is defined as \ref{Eq-def-tildef}.
\end{definition}
Note that 
$$
\|DF\|_{L^2((\lambda\times\nu\times\P))}^2=\sum_{n=1}^\infty n n!\|f_n\|_{L^2((\lambda\times\nu)^n)}^2. 
$$
Thus, $F\in \mD^1_2$ if and only if $F\in L^2(\bP)$ and $DF\in L^2((\lambda\times \nu\times \bP))$. 

\medskip

The next lemma about the closability of Malliavin derivative is a consequence of \cite[Theorem 12.6]{di2009malliavin}.  
\begin{lemma}
Assume $F_n\in \mathbb{D}^1_2$, $F_n\rightarrow F$ in $L^2(\P)$ and
$$\sup_n\|DF_n\|_{L^2(\lambda\times\nu\times \P)}\leq M<\infty. $$
Then, $F\in \mathbb{D}^1_2$ and 
$$\|DF\|_{L^2((\lambda\times\nu\times\P))}\leq M. $$
\end{lemma}\label{Le-Close}
\begin{proof}
By our assumption, $\{DF_n\}_{n\in \mN}$ is bounded in $L^2(\lambda \times \nu \times \bP)$, thus Banach-Saks theorem implies that, the Ces\`aro mean sequence of a suitable subsequence of $\{DF_n\}$, say $\{DF_{n_k}\}$,  converges strongly to some $G\in L^2(\lambda \times \nu \times \bP)$, i.e. 
$$
D \left(\frac{1}{m}\sum_{k=1}^m F_{n_k}\right) =\frac{1}{m}\sum_{k=1}^m DF_{n_k} \to G, \ \mbox{in}\ L^2(\lambda \times \nu \times \bP). 
$$
On the other hand, $\frac{1}{m}\sum_{k=1}^m F_{n_k} \to F$ in $L^2(\bP)$, by \cite[Theorem 12.6]{di2009malliavin}, we get $F\in \mD^1_2$, $DF=G$ and 
\begin{align*}
&\|DF\|_{L^2((\lambda\times\nu\times\P))} =\lim_{m\to\infty} \frac{1}{m} \left\|\sum_{k=1}^m DF_{n_k}\right\|_{L^2((\lambda\times\nu\times\P))} \\
 \leq& \liminf_{m\to\infty} \frac{1}{m}  \sum_{k=1}^m  \|DF_{n_k}\|_{L^2((\lambda\times\nu\times\P))} \leq M. 
\end{align*}
\end{proof}

\section{A study of nonlocal parabolic equations}

In this section we study the solvability and regularity of nonlocal elliptic equations with first order terms. First of all, we introduce the nonlocal operator studied in this work. 
Given a Borel measurable function $\sigma:  \mR^d\to\mR^d\otimes\mR^d$, we define 
$$
\cL f(x)=L_\sigma f(x):= \int_{\mR^d}\Big(f(x+\sigma(x) z)-f(x)-\nabla f(x)\cdot \sigma(x) z\1_{B_1}(z) \Big)\nu(\dif z) . 
$$
When $\si(x)=\si$ is a matrix, by Fourier's transform, we have
$$
\widehat{L_\sigma f}(\xi)=-\psi_\sigma(\xi)\hat f(\xi),
$$
where the symbol $\psi_\sigma(\xi)$ takes the form
$$
\psi_\sigma(\xi)=-\int_{\mR^d}(\e^{\mathrm{i}\xi\cdot \sigma z}-1-\mathrm{i}\sigma z\1_{B_1}(z) \cdot\xi)\nu(\dif z).
$$
Next we consider the following equation: 
\begin{align}\label{Eq-Resolvent}
\lambda u-\cL u -b\cdot\nabla u=f,\ \ \lambda> 0. 
\end{align}

\subsection{Constant coefficient case: $\sigma(x)=\sigma$}
In this subsection we consider equation \eqref{Eq-Resolvent} with non-degenerate constant coefficient $\sigma(x)=\sigma\in \mR^{d\times d}$. First of all, we establish the 
following Bernstein type inequality for nonlocal operator $L_\sigma$, which will play a crucial role in the sequel.  Recall that $\Delta_j$ is the dyadic block operator defined in section 2. 

\bl\label{Le-FracBI}
Let $p\geq 2$ and  $j=-1, 0, 1, \cdots$. Assume that $\nu$ satisfies \eqref{Eq-nu-lower} and $f\in \sS'(\mR^d)$ with $\Delta_j f\in L^p$, then 
$$
\int_{\mR^d}|\Delta_{j} f|^{p-2} \Delta_{j} fL_\sigma \Delta_{j} f \,\d x\leq 0,
$$
Moreover, there are constants $c>0, j_0\in \mN$ only depending on $d, p, c_0,\rho$ and $\Lambda$ such that for any $j=j_0, j_0+1, j_0+2,\cdots$,
\begin{align}\label{Eq-Bernstein2}
-\int_{\mR^d}|\Delta_j f|^{p-2} \Delta_j fL_\sigma \Delta_j f \,\d x\geq c 2^{\alpha j}\|\Delta_j f\|_p^p. 
\end{align}
\el

\begin{proof}
Without loss of generality, we can assume $\si =\mI$ and $f, g\in \sS(\mR^d)$. Denote $L=L_{\mI}$ and $\<f,g\>:=\int_{\mR^d} f(x) g(x) \dif x$.  By the symmetry
of $\nu$,  we have 
\begin{align}\label{eq-symmetry}
\<f, -L g\> =\frac{1}{2} \int_{\mR^d\times\mR^d} (f(x+z)-f(x)) (g(x+z)-g(x)) \nu(\d z) \d x. 
\end{align}
By \eqref{eq-symmetry} and the following elementary inequality: 
$$p(a-b)(|a|^{p-2}a-|b|^{p-2}b)\geq (a|a|^{\frac{p}{2}-1}-b|b|^{\frac{p}{2}-1})^2,\,\,\, \forall p\geq2, \,a,b\in \mR, $$
we get, 
\begin{align}\label{Eq-33}
\begin{aligned}
&\< f|f|^{p-2}, -L f\>\\
=& \frac{1}{2}\int_{\mR^d\times\mR^d}(f(x+ z)-f(x))(|f|^{p-2}f(x+ z)-|f|^{p-2}f(x))\nu(\d z)\d x\\
\geq &\frac{1}{2p} \int_{\mR^d\times\mR^d} (f|f|^{\frac{p}{2}-1}(x+ z)-f|f|^{\frac{p}{2}-1}(x))^2\nu(\d z)\d x\\
=& \frac{1}{p}\left\<f|f|^{\frac{p}{2}-1}, -L (f|f|^{\frac{p}{2}-1})\right\>=\frac{1}{p} \int_{\mR^d} \Big|(-L)^{\frac{1}{2}} (f|f|^{\frac{p}{2}-1})\Big|^2\d x.
\end{aligned}
\end{align}
This implies our first desired inequality. Noting that for all $x\in[0,1]$, $1-\cos x\geq c x^2$, by \ref{A_nu}, for any $|\xi|>\rho^{-1}$, 
\begin{align*}
\psi (\xi)=&\int_{\mR^d}(1-\cos(z\cdot\xi))\nu(\d z)\geq c\int_{|z|\leq |\xi|^{-1}} | z\cdot \xi|^2 \nu(\d z) \\
\gtrsim & | \xi|^2 \int_{|z|\leq |\xi|^{-1}} \left|z\cdot \frac{\xi}{| \xi|}\right|^2 \nu(\d z) \overset{\eqref{Eq-nu-lower}}{\gtrsim}|\xi|^\a.  
\end{align*}
On the other hand, it is easy to see that 
$$
\psi (\xi)\gtrsim |\xi|^2, \quad \forall |\xi|\leq \rho^{-1}.  
$$
Thus, $\psi(\xi)\gtrsim |\xi|^\a-1$.  Noting that $\sF((-L)^{1/2} g)(\xi)=\sqrt{\psi(\xi)} \sF(g)(\xi)$, by  Plancherel formula, 
\begin{equation}\label{AEq1}
\begin{split}
&\int_{\mR^d} \Big|(-L)^{\frac{1}{2}} (f|f|^{\frac{p}{2}-1})\Big|^2\d x=\int_{\mR^d}\psi  (\xi) \big|\sF(f|f|^{\frac{p}{2}-1})(\xi)\big|^2  \d \xi\\
\geq & c\int_{\mR^d}|\xi|^\a\big|\sF(f|f|^{\frac{p}{2}-1})(\xi)\big|^2  \d \xi- C\int_{\mR^d} \big|\sF(f|f|^{\frac{p}{2}-1})(\xi)\big|^2 \d \xi\\
\geq& c \int_{\mR^d} \Big|(-\Delta)^{\frac{\a}{4}} (f|f|^{\frac{p}{2}-1})\Big|^2\d x-C\int_{\mR^d} |f|^p\d x. 
\end{split}
\end{equation}
Combing \eqref{Eq-33}, \eqref{AEq1} and using the elementary inequality: 
$$\big|a|a|^{\frac{p}{2}-1}-b|b|^{\frac{p}{2}-1}\big|^2\geq c_p|a-b|^p, \,\, \forall a,b\in \mR,\,p\geq 2, $$
we obtain 
\begin{align*}
\int_{\mR^d} f|f|^{p-2}(-L f)\d x \geq &c \int_{\mR^d} \Big|(-\Delta)^{\frac{\a}{4}} (f|f|^{\frac{p}{2}-1})\Big|^2\d x-C \int_{\mR^d} |f|^p\d x\\
\geq & c\int_{\mR^d\times\mR^d} \frac{|f|f|^{\frac{p}{2}-1}(x)- f|f|^{\frac{p}{2}-1}(y)|^2}{|x-y|^{d+\a}}\d x\d y-C\|f\|_p^p\\
\geq&   c \int_{\mR^d\times\mR^d} \frac{|f(x)-f(y)|^p}{|x-y|^{d+\frac{\a}{p}\cdot p}}\d x\d y -C\|f\|_p^p= c [f]_{\frac{\a}{p},p}^p-C \|f\|_p^p. 
\end{align*}
Noting that $\frac{\a}{p}\leq \frac{\a}{2}<1$, by Theorem 2.36 of \cite{bahouri2011fourier}, for any $j\geq 0$, 
\begin{align*}
[\Delta_j f]_{\frac{\a}{p},p}^p=\|f\|_{\dot{B}^{\frac{\a}{p}}_{p,p}}=\sum_{k=-\infty}^{\infty} 2^{\a k} \|\Delta_k \Delta_j f\|_p^p\asymp 2^{\a j} \|\Delta_j f\|_p^p .\end{align*}
Thus, 
\begin{align*}
-\int_{\mR^d}|\Delta_j f|^{p-2} \Delta_j fL  \Delta_j f\dif x\geq (c 2^{\alpha j}-C)\|\Delta_j f\|_p^p. 
\end{align*}
Letting $j_0= 1+ \log_2 (C/c)/\a$, we get the desired result. 
\end{proof}
\br
Similar result was also proved in \cite{chen2021supercritical} by using inequality (1.5) in \cite{chen2007new}.  Our proof presented here is much simpler and self-contained. 
\er

Then we present a simple commutator estimate, which is a special case of Lemma 2.3 in \cite{chen2021supercritical}. 
\bl\label{Le-Comm}
Let $j\geq -1$, $\beta\in(0,1)$. Assume that $b\in C_b^\beta$ and $u\in W^1_p$, then 
\be\label{Eq-Comm}
\|[\Delta_j, b\cdot \nabla]u\|_p \lesssim 2^{-\beta j} \|b\|_{C_b^\beta} \|\nabla u\|_{p}. 
\ee
\el
\begin{proof}
By \eqref{Eq-conv-delta}, we have
\begin{align*}
[\Delta_j, b\cdot \nabla]u(x)=\int_{\mR^d}h_j(y)(b(x-y)-b(x))\cdot \nabla u(x-y)\dif y.  
\end{align*}
For $j\geq 0$, by Minkowski 's inequality and \eqref{Eq-hj}, we have
\begin{align*}
	\begin{split}
	\|[\Delta_j, b\cdot\nabla]u\|_p&\leq \int_{\mR^d}h_j(y)\|b(\cdot-y)-b(\cdot)\|_{\infty}\|\nabla u\|_{p}\dif y\\
	&\lesssim\|b\|_{C_b^\beta}\|\nabla u\|_{p}\int_{\mR^d}|h_j(y)|\,|y|^{\beta}\dif y\\
	&=\|b\|_{C_b^\beta}\|\nabla u\|_{p}  2^{-j\beta} \int_{\mR^d}|2h(2y)-h(y)|\,|y|^{\beta}\dif y\\
	&\lesssim 2^{-\beta j} \|b\|_{C_b^\beta}\|\nabla u\|_{p}. 
	\end{split}
\end{align*}
The estimate for $j=-1$ is simpler, so we complete our proof. 
\end{proof}

Last we state our main result of this subsection.
\bl\label{Priori1}
Let $\alpha\in(0,2)$,   $(1-\a)^+<\gamma <\beta<1$ and $p\in [2,\infty)$.
Suppose  that $\nu$ satisfies \eqref{Eq-nu-lower}, then there is a constant $\lambda_0>0$ such that for all $\lambda\geq \lambda_0$ and $f\in W^\gamma_{p}$, equation \eqref{Eq-Resolvent} admits a unique solution $u\in W^{\alpha+\gamma}_p$ . 
Moreover, 
\begin{equation}\label{Eq-Apriori1}
\lambda\|u\|_{W^\gamma_p}+\|u\|_{W^{\alpha+\gamma}_{p}}\leq C \|f\|_{W^{\gamma}_{p}},
\end{equation}
where $\lambda_0, C$ depend only on $d,p,\alpha,\beta,\gamma, c_0, \rho, \Lambda$ and $\|b\|_{C_b^{\beta}}$. 
\el
\begin{proof}
Thanks to method of continuity, our main task is to show the apriori estimates \eqref{Eq-Apriori1}. Assume $u\in W^{\a+\gamma}_p$ satisfies equation \eqref{Eq-Resolvent}. Using operator $\Delta_j$ act on both sides of \eqref{Eq-Resolvent} and noticing that $\Delta_j\cL=\Delta_j L_\si=L_\si \Delta_j =\cL \Delta_j$, we have
$$
\lambda \Delta_ju=\cL \Delta_j u+\Delta_j(b\cdot\nabla u)+\Delta_jf.
$$
For $p\geq 2$, by the chain rule or multiplying both sides by $|\Delta_j u|^{p-2}\Delta_j u$ and then integrating in $x$, we obtain
\begin{align*}
\lambda \int_{\mR^d} |\Delta_j u|^p=&\int_{\mR^d} |\Delta_j u|^{p-2}\Delta_j u\Big[\cL\Delta_j u+\Delta_j(b\cdot\nabla u)+\Delta_jf \Big]\dif x\\
=&\int_{\mR^d}|\Delta_j u|^{p-2}\Delta_j u\cL\Delta_j u\dif x+\int_{\mR^d}|\Delta_j u|^{p-2}\Delta_j u\,[\Delta_j,b\cdot\nabla]u\dif x\\
&+\!\!\int_{\mR^d}|\Delta_j u|^{p-2}\Delta_j u\,(b\cdot\nabla) \Delta_j u\dif x+\!\!\int_{\mR^d}|\Delta_j u|^{p-2}\Delta_j u\Delta_jf\dif x\\
=&: I^{(1)}_j+I^{(2)}_j+I^{(3)}_j+I^{(4)}_j.
\end{align*}
For $I^{(1)}_j$, recalling $\cL=L_{\sigma}$ and by Lemma \ref{Le-FracBI}, there are constants $c>0$ and $j_0\in \mN$ such that
\begin{align*}
I^{(1)}_{j}\leq 0,\ \ j\geq -1; \ \  I^{(1)}_j\leq -c2^{\alpha j}\|\Delta_j u\|^p_p,\ \ j=j_0,j_0+1,j_0+2,\cdots.
\end{align*}
For $I^{(2)}_j$, using Lemma \ref{Le-Comm} and  H\"older's inequality, we have that for all $j=-1,0,1,\cdots$,
\begin{align*}
I^{(2)}_j&\leq\|[\Delta_j,b\cdot\nabla]u\|_p\|\Delta_j u\|_p^{p-1}\\
&\lesssim 2^{-\beta j}\|b\|_{C_b^\beta}\|\nabla u\|_{p}
\|\Delta_j u\|_p^{p-1}.
\end{align*}
For $I^{(3)}_j$, we write
\begin{align*}
I^{(3)}_j&=\int_{\mR^d}((b-S_jb)\cdot\nabla) \Delta_j u\,|\Delta_j u|^{p-2}\Delta_j u\dif x\\
&+\int_{\mR^d}(S_jb\cdot\nabla) \Delta_j u\,|\Delta_j u|^{p-2}\Delta_j u\dif x=:I^{(31)}_j+I^{(32)}_j.
\end{align*}
For $I^{(31)}_j$, by Bernstein's inequality \eqref{Eq-BI1}, we have
\begin{align*}
I^{(31)}_j&\leq\sum_{k\geq j}\|(\Delta_kb\cdot\nabla) \Delta_j u\|_p\|\Delta_j u\|^{p-1}_p\\
&\leq\sum_{k\geq j} \|\Delta_kb\|_{\infty}\|\nabla\Delta_j u\|_{p} \|\Delta_j u\|^{p-1}_p\\
&\lesssim 2^j\|\Delta_j u\|^{p}_p\sum_{k\geq j} \|\Delta_kb\|_{\infty}\leq 2^j\|\Delta_j u\|^{p}_p \|b\|_{C_b^\beta}\sum_{k\geq j} 2^{-\beta k}\\
&\lesssim 2^{(1-\beta)j} \|b\|_{C_b^\beta}\|\Delta_j u\|^{p}_p.
\end{align*}
For $I^{(32)}_j$, by integration by parts formula and \eqref{Eq-BI1} again, we have
\begin{align*}
I^{(32)}_j&=\frac{1}{p}\int_{\mR^d}(S_jb\cdot\nabla) |\Delta_j u|^p\dif x=-\frac{1}{p}\int_{\mR^d}S_j\div b\, |\Delta_j u|^p\dif x\\
&\leq \frac{1}{p}\|S_j\div b\|_{\infty}\|\Delta_j u\|^p_{p}\leq \frac{1}{p}\sum_{k\leq j}\|\Delta_k \div b\|_{\infty}\|\Delta_j u\|^p_{p}\\
&\lesssim \sum_{k\leq j}2^{k} \|\Delta_kb\|_{\infty} \|\Delta_j u\|^p_{p}\\
&\lesssim 2^{(1-\beta)j} \|b\|_{C_b^\beta}\|\Delta_j u\|^p_{p}.
\end{align*}
Combining the above calculations, we obtain
\begin{align*}
 \lambda \|\Delta_j u\|_p^p +c2^{\alpha j}\1_{\{j\geq j_0\}}\|\Delta_j u\|^p_p \leq &C 2^{-\beta j}\|b\|_{C_b^\beta}\|\nabla u\|_{p}
\|\Delta_j u\|_p^{p-1}\\
&+C2^{(1-\beta)j} \|b\|_{C_b^\beta}\|\Delta_j u\|^p_{p}+C\|\Delta_j u\|^{p-1}_{p}\|\Delta_j f\|_p. 
\end{align*}
By dividing both sides by $\|\Delta_j u\|_p^{p-1}$, we get
\begin{align*}
&\lambda\|\Delta_j u\|_p+c2^{\alpha j}\1_{\{j\geq j_0\}}\|\Delta_j u\|_p-C2^{(1-\gamma)j} \|b\|_{C_b^\beta}\|\Delta_j u\|_{p}\\
\leq & C2^{-\beta j}\|b\|_{C_b^\beta}\|\nabla u\|_{p}+ C\|\Delta_j f\|_p. 
\end{align*}
Since $1-\gamma<\alpha$, for some $\lambda$ sufficiently large and all $j\geq -1$,
\begin{align}\label{Eq-KeyE}
{\lambda}\|\Delta_j u\|_p+ 2^{\alpha j} \|\Delta_j u\|_p\leq C2^{-\beta j} \|\nabla u\|_{p} + C\|\Delta_j f\|_p. 
\end{align}
Multiplying both sides of \eqref{Eq-KeyE} by $2^{\gamma j}$ and then taking $\ell^p$ norm over $j$, we obtain
\begin{align*}
\lambda \|u\|_{W^\gamma_p}+\|u\|_{W^{\a+\gamma}_p} \leq C_1\Big(\|\nabla u\|_{p}+\|f\|_{W^\gamma_p}\Big), 
\end{align*}
where $C_1$ only depends on $d,p, \a,\beta,\gamma, c_0, \rho$ and  $\|b\|_{C_b^\beta}$. Recalling that $\a+\gamma>1$ and using interpolation theorem, we have 
$\|\nabla u\|_{p}\leq \frac{1}{2C_1} \|u\|_{W^{\a+\gamma}_p}+C' \|u\|_{W^\gamma_p}$. 
Choosing $\lambda_0>2C_1C'$,  we complete the proof for \eqref{Eq-Apriori1}. 
\end{proof}

\subsection{Varying coefficient case} In this subsection, we consider the varying coefficient case.   
\bl
Assume that $\a\in (0,2)$ and $\nu$ satisfies condition \eqref{nu-upper}, then for any  $r\in (0,\rho)$ and $0\leq \theta'<\a<\theta$, 
\begin{align}\label{Eq-thetaorder}
\int_{|z|\leq r} |z|^\theta \nu (\d z) \lesssim r^{\theta-\a}  
\end{align}
and
\begin{align}\label{Eq-thetaorder'}
\int_{r<|z|\leq 1} |z|^{\theta'} \nu (\d z) \lesssim r^{\theta'-\a}. 
\end{align}
\el
\begin{proof}
Noting $\theta>\a$, one sees that 
\begin{align*}
\int_{|z|\leq r} |z|^\theta \nu (\d z) \leq &\sum_{k=0}^\infty (2^{-k} r)^{\theta-2} \int_{2^{-k-1}r<|z|\leq 2^{-k} r} |z|^2\nu(\d z)\\
\lesssim & r^{\theta-\a} \sum_{k=0}^\infty 2^{k(\a-\theta)} \lesssim r^{\theta-\a}. 
\end{align*}
Similarly, 
\begin{align*}
\int_{r<|z|\leq 1} |z|^{\theta'} \nu (\d z)& \leq \int_{r<|z|\leq \rho} |z|^{\theta'} \nu (\d z) + \rho^{\theta'}\int_{\rho<|z|\leq 1}  \nu (\d z)\\
& \lesssim \sum_{k=0}^{\log_2(\rho/r)} (2^{-k} \rho)^{\theta'-2} \int_{2^{-k-1}\rho<|z|\leq 2^{-k} \rho} |z|^2\nu(\d z)+1\lesssim r^{\theta'-\a}. 
\end{align*}
\end{proof}

We also need the following lemma (see \cite[Theorem 2.4.7]{triebel1992theory}) in order to localize the resolvent equation. 

\bl[localization principle]\label{Le-Partition}
Let $\zeta_k\in C_c^\infty, k=1,2,\cdots$. Assume for any multi-index $\a$, \\
$\sup_{x\in \mR^d} \sum_k |\p^\a\zeta_k(x)|\leq C_\a<\infty$. Then there is a constant $C$ such that 
$$
\sum_{k} \|u\zeta_k\|_{W^s_p}^p \leq C \|u\|_{W^{s}_p}^p. 
$$
Moreover, if $\sup_{x\in \mR^d}\sum_{k} |\zeta_k(x)|^p\geq c>0$, 
then  we have 
\begin{align}
\|u\|^p_{W^s_p}\asymp \sum_k \|u\zeta_k\|^p_{W^s_p}. 
\end{align}
\el
The following lemma is taken from \cite[Lemma 3.5]{mikulevivcius1992cauchy}.
\bl\label{Le35}Suppose $s\in(0,2)$, $p>\max\{1, \tfrac{d}{s}\}$, then
\begin{align}\label{Eq-MaxLp}
\left\|\sup_{y\neq 0}\frac{|f(\cdot+y)-f(\cdot)-\nabla f(\cdot) \cdot y^{(s)}|}{|y|^s}\right\|_p\leq C \|f\|_{H^s_p}, 
\end{align}
where $y^{(s)}= y\1_{s\in (1,2)} + y\1_{s=1}\1_{B_1}(y)$. 
\el

\bt\label{Th-Lp-Est}
Let $\alpha\in(0,2)$,  $(1-\a)^+<\gamma<\beta< 1$ and $p> \frac{d}{\a\wedge 1}\vee\frac{d}{\a+\gamma-1}$. Assume $\nu$ is supported on $B_R$ satisfying conditions \ref{A_nu}, $\si$ satisfies \ref{A_si} and $b\in C_b^\beta$. Then there is a constant $\lambda_0$  such that for any $\lambda\geq \lambda_0$ and $f\in \sW^\gamma_{p}$, equation \eqref{Eq-Resolvent} admits a unique solution in $\sW^{\a+\gamma}_p$. Moreover, 
\begin{align*}
\lambda\|u\|_{\sW^\gamma_p}+\|u\|_{\sW^{\alpha+\gamma}_{p}}\leq C \|f\|_{\sW^{\gamma}_{p}}, 
\end{align*}
where $\lambda_0, C$ depend only on $d,p,\alpha,\beta,\gamma, R, c_0, \rho, \Lambda $ and $\|b\|_{C_b^{\beta}}$.   
\et
When $\si$ is independent of $x$, then we do not need to assume $\nu$ is compactly supported (see Remark \ref{Rek-Const}). However, if $\sigma$ does depend on $x$ and 
the L\'evy measure $\nu$ is not necessarily absolutely continuous w.r.t. the Lebesgue measure, it seems not feasible to  show that for any $f\in W^\gamma_{p}$,
$$
x\mapsto \int_{|z|\geq R}f(x+\sigma(x) z)\nu(\dif z)\in W^\gamma_{p}, 
$$
which is very essential if one wants to use the perturbation argument. So in this case, we need to assume $\mathrm{supp}\,\nu\subseteq B_R$ with some constant $R\geq 1$.  Thus, $\cL f(x)$ can be rewritten as   
\begin{align}\label{TrL}
\begin{aligned}
\cL f(x)= \int_{B_R}\Big(f(x+\sigma(x) z)-f(x)-\nabla f(x)\cdot \sigma(x) z\1_{B_1}(z) \Big)\nu(\dif z). 
\end{aligned}
\end{align}

In order to prove Theorem \ref{Th-Lp-Est}, we need a commutator estimate under the following assumption: there are  $\eps\in(0, 1)$ and $\Lambda\geq 1$ such that
\begin{align}\label{Eq-si-osc}
\begin{aligned}
	|\sigma(x)-\sigma(y)| \leq \Lambda|x-y|, \sigma(x)=\sigma(0),\ \ \forall |x|\geq\eps,\\
	\Lambda^{-1} |\xi|^2\leq|\sigma(0)\xi|^2\leq\Lambda|\xi|^2,\ \forall \xi\in\mR^d.
\end{aligned}
\end{align}
If ${\rm supp}\, \nu\subseteq B_R$, we denote 
$$
L_{\sigma}^\a f(x) :=\int_{B_R}\Big(f(x+\sigma(x) z)-f(x)-\nabla f(x)\cdot \sigma(x) z\1_{\alpha\in [1,2)}\Big)\nu(\dif z). 
$$
\bl\label{Le38}
Suppose $\nu$ is supported on $B_R$ satisfying \ref{A_nu} and $\sigma$ satisfies \eqref{Eq-si-osc}.  Then for any  $s\in(0,1), p>1$, we have
\begin{align}\label{Eq-Com-LR}
\left\| \left[\Delta^{s/2},L^\a_\sigma\right]u\right\|_p\leq C \eps^{1-s+\frac{d}{p}} \left\{ 
\begin{aligned}
 &\| u\|_{C_b^1}, \ &\a\in (0,1);\\
 &\| u\|_{C_b^{\a+\delta}}, \ & \a\in[1,2), \, \delta>0, 
\end{aligned}
\right. 
\end{align}
where $[\Delta^{s/2},L^\a_\sigma]u:=\Delta^{s/2}L^\a_\sigma u-L^\a_\sigma \Delta^{s/2}u$, and the constant $C>0$ is independent of $\eps$.
\el
The reader can find the proof of above lemma  in \cite{chen2021supercritical}.

\bl\label{Comm2}
Suppose  $\nu$ is supported on $B_R$ satisfying \ref{A_nu} and $\sigma$ satisfies \eqref{Eq-si-osc}. Then for any $ \gamma\in(0,1)$,  $p\in(\tfrac{d}{\alpha\wedge 1},\infty)$,we have
\begin{align*}
\big\|({L^\a_{\sigma}}-{L^\a_{\sigma(0)}})f\big\|_{W^{\gamma}_{p}}
\leq c_\eps \|f\|_{W^{\alpha+\gamma}_p},
\end{align*}
where $c_\eps\to 0$ as $\eps\to 0$.
\el
\begin{proof}
We only prove the estimate for $\alpha\in(0,1)$. The case $\alpha\in[1,2)$ is similar. Without loss of generality, we can assume $R=1$. Define 	
$$
\cT_\sigma f:={L^\a_{\sigma}}f-{L^\a_{\sigma(0)}}f=\int_{|z|\leq 1}\Big(f(x+\sigma(x) z)-f(x+\sigma(0) z)\Big)\nu(\dif z).
$$
For $\theta\in(d/p,1]$, $\delta>0$ and $y_0\in\mR^d$, by \eqref{Eq-MaxLp}, we have
$$
\left\|\sup_{|y-y_0|\leq \delta}|f(\cdot+y)-f(\cdot+y_0)|\right\|_p\lesssim \delta^\theta\|f\|_{H^\theta_p},
$$
which implies that
\begin{align*}
&\|f(\cdot+\sigma(\cdot) z)-f(\cdot+\sigma(0) z)\|_p\no\\
&\leq\left\|\sup_{|y-\sigma(0) z|\leq\Lambda\eps |z|}|f(\cdot+y)-f(\cdot+\sigma(0) z)\right\|_p
\lesssim\eps^{\theta}|z|^\theta\|f\|_{H^\theta_p}.
\end{align*}
Thus, due to $p>\frac{d}{\alpha}$, we can choose $\frac{d}{p}<\theta'<\alpha<\theta\leq 1$ such that
\begin{align*}
\|\cT_\sigma f\|_p\lesssim \int_{|z|\leq 1}(\eps^{\theta}|z|^\theta\|f\|_{H^\theta_p})\wedge (\eps^{\theta'}|z|^{\theta'}\|f\|_{H^{\theta'}_p})\nu(\dif z), 
\end{align*}
which together with Bernstein's inequality  implies 
\begin{equation}\label{eq-Tf}
\begin{aligned}
\|\cT_\sigma \Delta_j f\|_p\lesssim &\int_{|z|\leq 1}(\eps^{\theta}|z|^\theta \|\Delta_j f\|_{H^\theta_p})\wedge (\eps^{\theta'}|z|^{\theta'}\|\Delta_j f\|_{H^{\theta'}_p})\nu(\dif z)\\
\lesssim&
\int_{|z|\leq 1}(\eps^{\theta}|z|^\theta 2^{\theta j})\wedge (\eps^{\theta'}|z|^{\theta'}2^{\theta' j})\nu(\dif z)\ \|\Delta_j f\|_p. 
\end{aligned}
\end{equation}
By \eqref{Eq-thetaorder} and \eqref{Eq-thetaorder'}
\begin{align*}
&\int_{|z|\leq 1}(\eps^{\theta}|z|^\theta 2^{\theta j})\wedge (\eps^{\theta'}|z|^{\theta'}2^{\theta' j})\nu(\dif z)\\
=&\eps^{\theta}2^{\theta j}\int_{|z|\leq 2^{-j}}|z|^\theta\nu(\dif z)+\eps^{\theta'} 2^{\theta' j}\int_{2^{-j}<|z|\leq 1}|z|^{\theta'}\nu(\dif z)\lesssim\eps^{\theta'}2^{\alpha j}.
\end{align*}
Substituting the above estimate into \eqref{eq-Tf}, we obtain 
\be\label{eq-TDf}
\|\cT_\sigma \Delta_j f\|_p\lesssim \eps^{\theta'} 2^{\a j} \|\Delta_j f\|_{p}. 
\ee
Note that 
$$
\|\Delta_i\cT_\sigma f\|_p
\leq\sum_{j>i}\|\Delta_i\cT_\sigma\Delta_jf\|_p+\sum_{j\leq i}\|\Delta_i\cT_\sigma\Delta_jf\|_p=:\cJ_i^{(1)}+\cJ_i^{(2)}.
$$
For $\cJ_i^{(1)}$, by \eqref{eq-TDf}, we have
$$
\cJ_i^{(1)}\lesssim \eps^{\theta'}\sum_{j>i}2^{\alpha j}\|\Delta_j f\|_p\lesssim \eps^{\theta'} \|f\|_{W^{\alpha+\gamma}_{p}} \sum_{j>i}2^{-\gamma j}c_j, 
$$
where $c_j: =2^{(\a+\gamma)j}\|\Delta_jf\|_p/ \|f\|_{W^{\a+\gamma}_p}$. Hence, 
\begin{align*}
2^{\gamma i}\cJ^{(1)}_i\lesssim& \eps^{\theta'} \|f\|_{W^{\a+\gamma}_p} \sum_{j\in \mZ}   \1_{\{(i-j)<0\}} 2^{\gamma(i-j)} \cdot \1_{\{j\geq -1\}}c_j \\
=& \eps^{\theta'} \|f\|_{W^{\a+\gamma}_p} (a*b)_i, 
\end{align*}
where $a_k=\1_{\{k\geq -1\}} c_k$ and $b_k=\1_{\{k<0\}}2^{\gamma k}$ $(\forall k\in \mZ)$. This yields  
\begin{align}\label{J1}
\begin{aligned}
& \|2^{\gamma i} \cJ^{(1)}_i\|_{\ell^p}\lesssim \eps^{\theta'} \|f\|_{W^{\a+\gamma}_p}\|a\|_{\ell^p} \|b\|_{\ell^1}\\
 \lesssim&\eps^{\theta'} \|f\|_{W^{\a+\gamma}_p} \Big(\sum_{j\geq -1} c_j^p\Big)^{1/p} \lesssim \eps^{\theta'} \|f\|_{W^{\a+\gamma}_p}. 
\end{aligned}
\end{align}

For $\cJ^{(2)}_i$, if $i=-1$, then 
$$
\cJ^{(2)}_{-1} = \|\Delta_{-1} \cT_\si \Delta_{-1}f\|_p\lesssim \eps^\a \|\Delta_{-1}f\|_{H^\a_p} \lesssim \eps^\a \|f\|_{W^{\a+\gamma}_p} c_{-1}
$$
If $i\geq 0$, choose $s\in(\gamma\vee(1-\a+d/p),1)$ in Lemma \ref{Le38}. By Bernstein's inequality and Lemma \ref{Le38}, we have 
	\begin{equation*}
	\begin{aligned}
	\cJ^{(2)}_i&= \sum_{-1\leq j\leq i}\|\Delta_i\Delta^{-s/2}\Delta^{s/2}\cT_\sigma\Delta_jf\|_p
	\overset{\eqref{Eq-BI2}}{\lesssim} 2^{-s i} \sum_{-1\leq j\leq i}\|\Delta^{s/2}\cT_\sigma\Delta_jf\|_p\\
	&{\leq} 2^{-s i} \sum_{-1\leq j\leq i}\left(\|[\Delta^{s/2},\cT_\sigma]\Delta_jf\|_p+\|\cT_\sigma\Delta^{s/2}\Delta_jf\|_p\right)\\
	&= 2^{-s i} \sum_{-1\leq j\leq i}\left(\|[\Delta^{s/2}, L^\a_\si ]\Delta_jf\|_p+\|\cT_\sigma\Delta_j \Delta^{s/2}f\|_p\right)\\
	&\overset{\eqref{Eq-Com-LR}, \eqref{eq-TDf}}{\lesssim }2^{-s i} c_\eps \sum_{-1\leq j\leq i}\left(\|\Delta_jf\|_{C_b^1}
	+ \|\Delta^{s/2}\Delta_jf\|_{H^\alpha_p}\right)\\
	&\overset{\eqref{Eq-BI1}}{\lesssim} 2^{-s i} c_\eps \sum_{-1\leq j\leq i}\left(2^{(1+d/p)j}\|\Delta_jf\|_p+ 2^{(s+\alpha)j} \|\Delta_jf\|_p\right)\\
	&\lesssim c_\eps\sum_{-1\leq j\leq i} 2^{(\alpha+s)j}\|\Delta_jf\|_p\\
	&\lesssim  c_\eps\|f\|_{W^{\a+\gamma}_p} 2^{-s i} \sum_{-1\leq j\leq i} 2^{(s-\gamma)j}c_j
	\end{aligned}
\end{equation*}
	Denoting $a_k=\1_{\{k\geq -1\}} c_k$,  $d_k= \1_{\{k\geq 0\}}2^{(\gamma-s)k}$, then $$
2^{\gamma i}\cJ_i^{(2)}\leq c_\eps\|f\|_{W^{\a+\gamma}_p} \sum_{-1\leq j\leq i}2^{(\gamma-s)(i-j)} c_j \leq c_\eps\|f\|_{W^{\a+\gamma}_p}  (d* a)_i. 
$$
Thus, 
\be\label{J2}
\|2^{\gamma i}\cJ_i^{(2)}\|_{\ell^p} \lesssim c_\eps\|f\|_{W^{\a+\gamma}_p} \|d\|_{\ell^1} \|a\|_{\ell^p}\lesssim c_\eps\|f\|_{W^{\a+\gamma}_p}. 
\ee
By \eqref{J1} and \eqref{J2}, we get 
\begin{align*}
\|\cT_\si f\|_{W^{\gamma}_p} \asymp& \|2^{\gamma i} \|\Delta_i \cT_\si f\|_p\|_{\ell^p} \leq \| 2^{\gamma i}\cJ_i^{(1)}\|_{\ell^p}+\|2^{\gamma i}\cJ_i^{(2)}\|_{\ell^p} \\
\lesssim & c_\eps\|f\|_{W^{\a+\gamma}_p}.  
\end{align*}
This yields our desired result. 
\end{proof}

Now we are on the position to prove Theorem  \ref{Th-Lp-Est}. 
\begin{proof}[Proof of Theorem  \ref{Th-Lp-Est}]
Like before, we only give the a prior estimate here. Assume $\{\zeta_k\}_{k\in \N}$ is a standard partition of unity and for each $k$, the support of $\zeta_k$ lies in a ball $B_{\eps/8}(y_k)$, and $\eps>0$ will be determined later. Also for any $k$, we take functions $\eta_k,\xi_k \in C^\infty_c(\mR^d)$ such that $\eta_k=1$ on $B_{\eps/4}(y_k)$, $\eta_k=0$ outside  $B_{\eps/2}(y_k)$, and $0\leq \eta_k\leq 1$; $\xi_k=1$ on  $B_{\eps/2}(y_k)$, $\xi_k=0$ outside  $B_\eps(y_k)$, and $0\leq \xi_k\leq 1$. Define $\si_k(x)=\xi_k(x)\si(x)+(1-\xi_k(x))\si(y_k)$,
$$
\cL_kf(x):= \int_{B_R } \Big(f(x+\si_k(x)z)-f(x)-\nabla f(x)\sigma_k(x) z\1_{B_1}(z)\Big)\nu(\d z). 
$$
$$
L_{k}f(x):=\int_{B_R } \Big(f(x+\si(y_k)z)-f(x)-\nabla f(x)\sigma(y_k) z\1_{B_1}(z)\Big)\nu(\d z). 
$$
Multiplying $\zeta_k$ on both side of \eqref{Eq-Resolvent},  we get 
\begin{equation}
\lambda( u\zeta_k)-L_k (u\zeta_k)- b\cdot\nabla (u\zeta_k)=f\zeta_k+\zeta_k (b\cdot\nabla u)-b\cdot\nabla (u\zeta_k)+\zeta_k (\cL u)-L_{k}(\zeta_ku). 
\end{equation}
Lemma \ref{Priori1} yields, 
$$
\lambda \|u\zeta_k\|_{W^{\gamma}_p}+\|u\zeta_k\|_{W^{\a+\gamma}_p}\lesssim \left(\|f\zeta_k\|_{W^\gamma_p}+\|ub\cdot \nabla \zeta_k\|_{W^\gamma_p}+\|\zeta_k (\cL u)-L_{k}(\zeta_ku)\|_{W^\gamma_p}\right). 
$$
So by Lemma \ref{Le-Partition}, we have, 
\be\label{Sum}
\begin{split}
\lambda^p\|u\|_{W^\gamma_p}^p+\|u\|_{W^{\a+\gamma}_p}^p\lesssim & \sum_k \lambda^p\|u\zeta_k\|_{W^\gamma_p}^p+\|u\zeta_k\|_{W^{\a+\gamma}_p}^p\\
\lesssim & \sum_k \left(\|f\zeta_k\|_{W^\gamma_p}^p+\|ub\cdot \nabla \zeta_k\|_{W^\gamma_p}^p+\|\zeta_k (\cL u)-L_{k}(\zeta_ku)\|_{W^\gamma_p}^p\right)
\end{split}
\ee
Again by Lemma \ref{Le-Partition}, \begin{align}\label{Sum-f}
\sum_k \|f\zeta_k\|^p_{W^\gamma_p}\asymp \|f\|_{W^\gamma_p}^p, \,\,\, \sum_k\|ub\cdot \nabla \zeta_k\|^p_{W^\gamma_p}\lesssim \|ub\|_{W^\gamma_p}^p\lesssim \|u\|_{W^\gamma_p}^p, 
\end{align} 
the last inequality above is due to the following fact: 
\begin{align*}
 \|ub\|_{W^\gamma_p}^p \lesssim&  \|ub\|_p^p+  \iint_{\mR^d\times\mR^d} \frac{|ub(x)-ub(y)|^p}{|x-y|^{d+\gamma p}} \d x\d y\\
 \leq& \|u\|_p^p \|b\|_{L^\infty}^p \d x\d y + \iint_{\mR^d\times \mR^d} \frac{|u(x)-u(y)|^p|b(x)|^p}{|x-y|^{d+\gamma p}}\\
&+ \iint_{\mR^d\times \mR^d} \frac{|u(y)|^p|b(x)-b(y)|^p}{|x-y|^{d+\gamma p}} \d x\d y\\
\lesssim& \|u\|_{p}^p\|b\|_{L^\infty}^p+ \|b\|_{L^\infty}^p  [u]_{\gamma ,p}^p + \int_{\mR^d} |u(y)|^p \d y \int_{\mR^d} \frac{[b]_{\beta}^p|x-y|^{\beta p}\wedge \|b\|_{L^\infty}^p}{|x-y|^{d+\gamma p}}\d x\\
\lesssim & \|u\|_{W^\gamma_p}^p \|b\|_{C_b^\beta}^p. 
\end{align*}
Next we estimate the third term in the last line of \eqref{Sum}. We only give the proof for $\a<1$ below, since the proof for $\a\geq 1$ is almost the same. Rewrite 
\begin{align*}
& \zeta_k(x) (\cL u)(x)-L_{k}(\zeta_ku)(x)\\
=&\left[\cL(u\zeta_k)(x)-L_k(u\zeta_k)(x)\right]\eta_k(x)+\left[\cL(u\zeta_k)(x)-L_k(u\zeta_k)(x)\right](1-\eta_k(x))\\
&- \left\{ u(x)\cL\zeta_k(x)+\int_{B_R} [u(x+\si(x)z)-u(x)]\cdot [\zeta_k(x+\si(x)z)-\zeta_k(x)]\nu(\d z)\right\}\\
=&: I^{(1)}_k(x)+I^{(2)}_k(x)-I^{(3)}_k(x). 
\end{align*}
For $I^{(1)}_k$, notice $\si_k(x)=\si(x)$ when $x$ is in the support of $\eta_k$, so we have
\begin{align*}
I^{(1)}_k(x)=&[\cL_k(u\zeta_k)(x)-L_k(u\zeta_k)(x)]\eta_k(x)\\
=& [L_{\sigma_k}^\a(u\zeta_k)(x)-L_{\sigma(y_k)}^\a(u\zeta_k)(x)]\eta_k(x). 
\end{align*}
By Lemma \ref{Comm2}, we have 
\begin{align}\label{I1}
\|I_k^{(1)}\|_{W^\gamma_p}\leq c_\eps\|u\zeta_k\|_{W^{\a+\gamma}_p} \ \ (c_\eps\to 0\ \mbox{as}\  \eps\to0). 
\end{align}
For $I_k^{(2)}(x)$, since $1-\eta_k(x)=0$ if $|x-y_k|\leq \frac{\eps}{4}$ and $u\zeta_k(x)=0$ if $|x-y_k|>\frac{\eps}{8}$, we have
\begin{align*}
I_k^{(2)}(x)=&\int_{\frac{\eps}{8\Lambda}\leq |z|<R } [u\zeta_k(x+\si(x)z)-u\zeta_k(x+\si(y_k)z)](1-\eta_k(x))\nu(\d z). 
\end{align*}
Choosing $s\in(d/p, 1\wedge(\a+\gamma-1))$, by Minkowski inequality, Lemma \ref{Le35} and interpolation theorem, we have 
\begin{align*}
\|I_k^{(2)}\|_p\leq & \int_{\frac{\eps}{8\Lambda}\leq |z| < R}\|u\zeta_k(\cdot+\si(\cdot)z)-u\zeta_k(\cdot+\si(y_k)z) \|_p\nu(\d z)\\
\lesssim &\int_{\frac{\eps}{8\Lambda}\leq |z| < R}  \Big\|\sup_{|y|\leq  2\Lambda |z|} \left|u\zeta_k(\cdot+y)-u\zeta_k(\cdot)\Big| \right\|_p \nu(\d z)\\
\overset{\eqref{Eq-MaxLp}}{\lesssim} & \|u\zeta_k\|_{H^s_p}\int_{\frac{\eps}{8\Lambda}\leq |z|<R}|z|^s  \nu(\d z)\lesssim \eps^{s-\a} \|u\zeta_k\|_{H^s_p}\\
\leq& \eps\|u\zeta_k\|_{W^{\a+\gamma}_p}+C_\eps\|u\zeta_k\|_p. 
\end{align*}
Moreover, noticing that 
\begin{align*}
\nabla I_k^{(2)}(x)
=&\int_{\frac{\eps}{8\Lambda}\leq |z|<R } \Big\{[\nabla(u\zeta_k)(x+\si(x)z)(\mI+\nabla \si(x)z)-\nabla(u\zeta_k)(x+\si(y_k)z)]\\
&(1-\eta_k(x))
- [u\zeta_k(x+\si(x)z)-u\zeta_k(x+\si(y_k)z)]\nabla\eta_k(x)\Big\}\nu(\d z)\\
=&\int_{\frac{\eps}{8\Lambda}\leq |z|<R } \Big\{ \nabla(u\zeta_k)(x+\si(x)z)\cdot  \nabla \si(x)z\ (1-\eta_k(x))\\
&\quad +[\nabla(u\zeta_k)(x+\si(x)z)-\nabla(u\zeta_k)(x+\si(y_k)z)](1-\eta_k(x))\\
&\quad- [u\zeta_k(x+\si(x)z)-u\zeta_k(x+\si(y_k)z)]\nabla\eta_k(x)\Big\}\nu(\d z), 
\end{align*}
like before, we choose $s\in(d/p, 1\wedge(\a+\gamma-1))$, then 
\begin{align*}
\|\nabla I^{(2)}_k\|_p\lesssim&  \int_{\frac{\eps}{8\Lambda}\leq |z|<R} \Bigg\{ \Big[ \|\nabla(u\zeta_k)(\cdot+\si(\cdot)z)-\nabla(u\zeta_k)(\cdot)\|_p+ \|\nabla(u\zeta_k)\|_p\Big] \cdot \|\nabla \si\|_\infty|z|\\
&+\|\nabla(u\zeta_k)(\cdot+\si(\cdot)z)-\nabla(u\zeta_k)(\cdot+\si(y_k)z)\|_p\\
&+\|\nabla\eta_k\|_{\infty}\cdot \|u\zeta_k(\cdot+\si(\cdot)z)-u\zeta_k(\cdot+\si(y_k)z)\|_p \Bigg\}  \nu(\d z)\\
\overset{\eqref{Eq-MaxLp}}{\lesssim_\eps }& \int_{\frac{\eps}{8\Lambda}\leq |z| < R}  (|z|^{1+s}+|z|^s) \left\|\sup_{y\in \mR^d} \frac{|\nabla (u\zeta_k)(\cdot+y)-\nabla(u\zeta_k)(\cdot)|}{|y|^s} \right\|_p \nu(\d z)\\
&+ \int_{\frac{\eps}{8\Lambda}\leq |z| < R}  |z| \|\nabla (u\zeta_k)\|_p \nu (\d z)\\
\lesssim_\eps& \|u\zeta_k\|_{H^{1+s}_p} \leq\eps\|u\zeta_k\|_{W^{\a+\gamma}_p}+C_\eps\|u\zeta_k\|_p. 
\end{align*}
So 
\begin{align}\label{I2}
\|I^{(2)}_k\|_{W^\gamma_p}\leq\|I^{(2)}_k\|_{W^1_p}\leq   \eps\|u\zeta_k\|_{W^{\a+\gamma}_p}+C_\eps\|u\zeta_k\|_p.
\end{align}
For $I^{(3)}_k(x)$, we have 
$$
\sup_{x\in \mR^d}\sum_k|\zeta_k(x+\si(x)z)-\zeta_k(x)|^p\lesssim |z|^p, \quad \forall |z|<R 
$$
and 
$$
\sup_{x\in \mR^d}\sum_k |(\cL\zeta_k)(x)|^p\lesssim 1. 
$$
Hence, for any $s\in(0,\a)$, 
\begin{align*}
&\left(\sum_k\|I^{(3)}_k\|_p^p\right)^{1/p}\leq\left(\sum_k \int_{\mR^d} |u(x)|^p|\cL\zeta_k(x)|^p\d x \right)^{1/p}\\
&+\int_{|z|<R } \left\{ \int_{\mR^d} |u(x+\si(x)z)-u(x)|^p\sum_k |\zeta_k(x+\si(x)z)-\zeta_k(x)|^p \d x \right\}^{1/p} \nu(\d z) \\
\lesssim&  \|u\|_p+\int_{|z|<R } \|u\|_{H^s_p}|z|^{s+1} \nu(\d z)\leq  \eps\|u\|_{W^{\a+\gamma}_p}+C_\eps\|u\|_{p}. 
\end{align*}
Similarly, we have  
$$
\left(\sum_k\|\nabla I^{(3)}_k\|_p^p\right)^{1/p}\leq \eps\|u\|_{W^{\a+\gamma}_p}+C_\eps \|u\|_p. 
$$
So,
\begin{align}\label{I3}
\left(\sum_k \|I^{(3)}_k\|_{W^\gamma_p}^p\right)^{1/p}\leq\left(\sum_k \|I^{(3)}_k\|_{W^1_p}^p\right)^{1/p} \leq   \eps\|u\|_{W^{\a+\gamma}_p}+C_\eps \|u\|_p.
\end{align}
Now by Lemma \ref{Le-Partition}, \eqref{Sum}-\eqref{I3}, and choosing $\eps$ sufficiently small and $\lambda_0$ sufficiently large, we get 
\be\label{eq-W-est}
\lambda\|u\|_{W^\gamma_p}+\|u\|_{W^{\alpha+\gamma}_{p}}\leq C \|f\|_{W^{\gamma}_{p}}. 
\ee
Now multiplying both sides of \eqref{Eq-Resolvent} by $\chi_z$, we have 
\begin{align*}
\lambda(u\chi_z)- \cL (u\chi_z) - b\cdot\nabla (u\chi_z) = g_z, 
\end{align*}
where 
$$
 g_z:= f\chi_z  + \chi_z\cL u-\cL (u\chi_z) -u b\cdot\nabla \chi_z. 
$$
We omit the subscript $z$ below. By definition, 
\begin{align*}
[\chi\cL u-\cL (u\chi) ](x)= \int_{B_R} u(x+\sigma(x)z)(\chi(x+\sigma(x)z)-\chi(x)) \nu(\d z), 
\end{align*}
so again by Lemma \ref{Le35}, 
\begin{align*}
\|\chi\cL u-\cL (u\chi) \|_p\leq&  \|u\|_\infty \left\|\int_{|z|<R} [\chi(\cdot+\sigma(\cdot) z)-\chi(\cdot)] \nu(\d z) \right\|_p\\
\leq & C \|u\|_\infty. 
\end{align*}
Noting that 
\begin{align*}
\nabla [\chi\cL u-\cL (u\chi)](x)= &\int_{B_R} \nabla u(x+\sigma(x)z) (\mI+ \nabla \sigma(x)z) (\chi(x+\sigma(x)z)-\chi(x)) \nu(\d z) \\
&+ \int_{B_R} u(x+\sigma(x)z) (\nabla\chi(x+\sigma(x)z)(\mI+\nabla \sigma(x)z)-\nabla\chi(x)) \nu(\d z), 
\end{align*}
we have 
\begin{align*}
\|\nabla[\chi\cL u-\cL (u\chi) ]\|_p\leq & C \|\nabla u\|_\infty  \left\|\int_{B_R} [\chi(\cdot+\sigma(\cdot) z)-\chi(\cdot)] \nu(\d z) \right\|_p\\
&+ C \|u\|_\infty \left\|\int_{B_R} [\nabla \chi(\cdot+\sigma(\cdot) z)-\nabla \chi(\cdot)] \nu(\d z) \right\|_p\\
&+ C \|u\|_\infty \left\|\int_{B_R} \nabla \chi(\cdot+\sigma(\cdot) z)\cdot \nabla \sigma(x)z \nu(\d z) \right\|_p \\
\leq & C \|u\|_{C_b^1} . 
\end{align*}
Thus, 
\begin{align*}
 \|g\|_{W^\gamma_p}\leq& \|f\chi\|_{W^\gamma_p}+ C \|u\|_{C_b^1}+ C\|b\|_{C_b^\beta} \|u\nabla \chi\|_{W^\gamma_p}\\
 \leq& C (\|f\|_{\sW^\gamma_p}+ \|u\|_{C_b^1}+ \|u\|_{\sW^\gamma_p}). 
\end{align*}
By \eqref{eq-W-est}, we get 
\begin{align*}
\|u\|_{\sW^{\a+\gamma}_p}+\lambda \|u\|_{\sW^{\gamma}_p}=&\sup_{z\in\mR^d}(\|u\chi_z\|_{W^{\a+\gamma}_p} + \lambda \|u
\chi_z\|_{W^{\gamma}_p}) \\
\leq & C\sup_{z\in \mR^d}\|g_z\|_{W^\gamma_p} \leq  (\|f\|_{\sW^\gamma_p}+ \|u\|_{C_b^1}+ \|u\|_{\sW^\gamma_p}). 
\end{align*}
By Lemma \ref{Le-embed} and interpolation theorem, for any $\delta>0$ there is a constant $C_\delta>0$ such that 
$$
\|u\|_{C_b^1} \leq \delta \|u\|_{\sW^{\a+\gamma}}+ C_\delta \|u\|_{\sW^\gamma}, 
$$
so we complete our proof by first choosing $\delta$ small and then choosing $\lambda_0$ sufficiently large. 
\end{proof}

\begin{proof}[Proof of Theorem \ref{Th-Holder}] Our Theorem \ref{Th-Holder} is a consequence of Lemma \ref{Le-embed} and Theorem \ref{Th-Lp-Est}. 
\end{proof}

\section{Proof of Theorem \ref{Th-Main}}
In this section, we prove our probabilistic results by using our analytic results proved in the previous section and Zvonkin's transform. Let $N(\dif t,\dif z)$ be the Poisson random measure associated with $Z$, that is,
$$
N((0,t]\times E)=\sum_{s\leq t}1_{E}(\Delta Z_s),$$
where $E$ is any compact set of $\mR^d\backslash \{0\}$ and $\Delta Z_s:=Z_s-Z_{s-}$. 
The  intensity measure of $N$ is denoted by $\dif t\nu(\dif z)$. Let $\widetilde N(\d t, \d z)=N(\d t, \d z)- \d t \nu(\d z)$, then 
$$
Z_t=\int^t_0\!\!\!\int_{|z|\leq 1} z \widetilde N(\dif s,\dif z)+ \int^t_0\!\!\!\int_{|z|>1} zN(\d s, \d z). 
$$
Thus, SDE \eqref{Eq-SDE} can be rewritten as
\begin{equation}\label{SDE3}
X_t=X_0+\int_0^t b(X_s)\d s+\int_0^t\!\!\!\int_{|z|\leq 1}\si(X_{s-})z \widetilde N(\dif s,\dif z)+\int_0^t\!\!\!\int_{|z|> 1}\si(X_{s-})z N(\dif s,\dif z) . 
\end{equation}
Now we are on the point of giving the proof of our main results. 
\begin{proof}[Proof of Theorem \ref{Th-Main}] 
(i). For the well-posedness, one can assume $\nu$ compactly supported on $B_{R}$ i.e.  $\sup_{t\geq 0}|\Delta Z_{t}|<R$, otherwise, we can take $\tau_0:=0, \tau_{k}:=\inf\{t>\tau_{k-1}: |\Delta Z_t|\geq R\}$ for any $k\geq 1$, and solve the SDE step by step. 

Now let $u$ be the solution of equation: 
$$
\lambda u-\cL u-b\cdot\nabla u=b, \quad \lambda\geq \lambda_0, 
$$
where 
$$
\cL u (x)= \int_{B_R} \left( u(x+z)-u(x)-\nabla u(x)\cdot \sigma(x)z\1_{B_1}(z) \right)\nu(\d z). 
$$
By Theorem \ref{Th-Holder} and interpolation theorem, for any $\mu\in(\alpha/2,\alpha+\beta-1)$, we have $u\in C_b^{1+\mu}$ with $\|u\|_{C_b^{1+\mu}}=c(\lambda, \mu)$ and  $c(\lambda, \mu)\to0$ as $\lambda\to \infty$. Choose $\lambda$ sufficiently large so that $\|\nabla u\|_\infty\leq 1/2$, then $\phi: x\mapsto x+u(x)$ is a $C_b^{1+\mu}$-diffeomorphism. By \eqref{SDE3} and a generalized version of It\^o's formula(c.f. \cite{priola2012pathwise}), we get 
\begin{align*}
u(X_t)=&u(X_0)+\int_0^t [\cL u+b\cdot \nabla u](X_s)\d s\\
&+\int_0^t\!\!\! \int_{B_R} [u(X_{s-}+\si(X_{s-})z)-u(X_{s-})]\widetilde{N} (\d s,\d z).
\end{align*}
Define $Y_t:= \phi(X_t)$, then 
\begin{equation}
\begin{split}
\label{Eq-Y}
Y_t&=\phi(X_t)=\phi(X_0)+\int_0^t  \lambda u(X_s)\d s\\
&+\int_0^t\!\!\!\int_{B_R} [\phi(X_{s-}+\si(X_{s-})z)-\phi(X_{s-})]\widetilde{N} (\d s,\d z)\\
&=Y_0+\int_0^t a(Y_s)\d s+\int_0^t\!\!\!\int_{B_R} g(Y_{s-},z)\widetilde{N} (\d s,\d z),
\end{split}
\end{equation}
where 
$$
a(y):=\lambda u (\phi^{-1}(y)); 
$$
\begin{align*}
g(y,z):=&\phi(\phi^{-1}(y)+\sigma(\phi^{-1}(y))z)-y\\
=&u(\phi^{-1}(y)+\sigma(\phi^{-1}(y))z)-u(\phi^{-1}(y))+\sigma(\phi^{-1}(y))z.
\end{align*}
It has been showed in \cite{priola2012pathwise} and \cite{chen2018stochastic} that in order to get the well-posedness of  \eqref{SDE3},  we only need to show the well-posedness of \eqref{Eq-Y}. 
Elementary calculations yield, 
$$
\nabla a(y)=\lambda \nabla u(\phi^{-1}(y))\nabla \phi^{-1}(y); 
$$
\begin{align}\label{Eq-de-g}
\begin{aligned}
\nabla_y g(y,z)=&[\nabla u(\phi^{-1}(y)+\si(\phi^{-1}(y))z)-\nabla u(\phi^{-1}(y))] \nabla\phi^{-1}(y)\\
&+\nabla u(\phi^{-1}(y)+\si(\phi^{-1}(y))z)\ \nabla\si(\phi^{-1}(y)) \ \nabla\phi^{-1}(y)z\\
&+\nabla \si(\phi^{-1}(y))\ \nabla\phi^{-1}(y) z. 
\end{aligned}
\end{align}
Fix $\mu\in(\alpha/2,\alpha+\beta-1)$.  Noting that $u\in C_b^{1+\mu}$ with $\|u\|_{C_b^{1+\mu}}=c(\lambda, \mu)$, 
we have, 
\begin{align}\label{Eq-De-a}
\|a\|_{C_b^{1+\mu}}<\infty , \quad |g(y,z)|\leq C |z|
\end{align}
and
\begin{align}\label{Eq-De-g}
\begin{aligned}
\|\nabla_y g(\cdot,z)\|_\infty\leq& \|\nabla u\|_{C_b^\mu} \|\nabla \phi^{-1}\|_\infty(\|\si\|_\infty \cdot |z|)^\mu\\
&+ \|\nabla u\|_\infty\|\nabla \si\|_\infty\|\nabla \phi^{-1}\|_\infty|z|+\|\nabla \si\|_\infty\|\nabla \phi^{-1}\|_\infty|z|\\
\leq & c(\lambda, \mu) (1-c(\lambda, \mu))^{-1} \big(\|\si\|_\infty^\mu  |z|^\mu+ \|\nabla\si\|_\infty|z| \big)+\|\nabla \si\|_\infty |z|. 
\end{aligned}
\end{align}
Thanks to the estimates \eqref{Eq-De-a} and \eqref{Eq-De-g}, the proof for existence and uniqueness of solution to \eqref{Eq-Y} becomes quite standard.  Let 
$$Y^0_t=Y_0;\quad Y^{n+1}_t=Y_0+\int_0^t a(Y^n_s)\d s +\int^t_0\!\!\!\int_{B_R} g(Y^n_{s-},z)\widetilde{N}(\d s,\d z), $$
then by Doob's inequality and \eqref{Eq-De-a}, \eqref{Eq-De-g}, we have
\begin{align*}
\E \sup_{0\leq s\leq t}|Y^{n+1}_s-Y_s^{n}|^2 \leq & C\|\nabla a\|_\infty^2 \E \int_0^t |Y^{n}_s-Y^{n-1}_s|^2\d s\\
&+ C \E \int^t_0\!\!\!\int_{B_R} \left|g(Y^n_{s-},z)-g(Y^{n-1}_{s-},z)\right|^2 \nu(\d z)\d s \\
\leq & C \left(\|\nabla a\|_\infty^2+ \int_{B_R} \|\nabla_y g(\cdot,z)\|_\infty^2  \nu(\d z)\right) \E \int_0^t |Y^n_s-Y^{n-1}_s|^2 \d s \\
\leq&  C\left[ 1+\int_{B_R} (|z|^{2\mu}+|z|^2) \nu(\d z) \right] \E \int_0^t |Y^n_s-Y^{n-1}_s|^2 \d s\\
\leq & C \E \int_0^t |Y^n_s-Y^{n-1}_s|^2 \d s\leq C t \E \sup_{0\leq s\leq t} |Y_s^n-Y_s^{n-1}|^2. 
\end{align*}
Choosing $T$ sufficiently small such that $CT\leq \frac{1}{2}$, we get 
\be\label{Eq-lim-Ynm}
\lim_{n, m\rightarrow\infty}\E\sup_{0\leq t\leq T} |Y^n_t-Y^{m}_t|^2=0. 
\ee
Note that all the estimates above do not  depends on the initial data $Y_0$, so we obtain that  \eqref{Eq-lim-Ynm} holds for any $T>0$. The limit point $Y$ of $\{Y^n\}_n$ is a strong solution to \eqref{Eq-Y}.  The uniqueness for \eqref{Eq-Y} can be obtained by using Gronwall's inequality and similar estimates as above. 
 
\medskip

Next we show that for each $t$, $X_t$ is Malliavin differentiable. Noting that $\nabla (\phi^{-1})(y)= (\nabla \phi)^{-1}\circ \phi^{-1}(y)$ and $(\nabla \phi)^{-1} =\mathrm{adj}(\nabla \phi) (\det \nabla \phi)^{-1}\in C_b^\mu$ ($\mathrm{adj}(\nabla \phi)$ is the adjugate matrix of $\nabla \phi$), we get $\phi^{-1}\in C_b^{1+\mu}$. Since $X_t=\phi^{-1}(Y_t)$, by Theorem 12.8 of \cite{di2009malliavin}, we only need to show that $Y_t$ is Malliavin differentiable. By the closability of Malliavin derivative (Lemma \ref{Le-Close}), the desired result can be derived from following estimate: 
\begin{equation}\label{DerivateBound}
\sup_{n\in \mN; t\in [0,T]} \|D_{r,z} Y_t^n\|_{L^2(\lambda\times\nu\times\P)}<\infty. 
\end{equation}
Assume $Y^n_t$ above is Malliavin differentiable  for each $t$ and $\sup_{t\in [0,T]} \|D_{r,z} Y_t^n\|_{L^2(\lambda\times\nu\times\P)}<\infty$. By \eqref{Eq-De-a}, we have 
$$
a(Y^n_s)\in L^2(\bP), \quad a(Y_s^n+D_{r,z}Y^n_s)-a(Y^n_s) \in L^2(\lambda\times \nu\times \bP). 
$$
Thanks to Theorem 12.8 of \cite{di2009malliavin}, we obtain $a(Y_s^n)\in \mD^1_2$ and 
$$
D_{r,z}a(Y^n_s)= a(Y^n_s+D_{r,z}Y^n_s)-a(Y^n_s). 
$$
Similarly, by \eqref{Eq-De-a}, 
$$
\bE \int_0^T\int_{\mR^d}|g(Y_{s-}, \eta)|^2 \nu(\d \eta) \d s<\infty, 
$$
and by \eqref{Eq-De-g} and It\^o's isometry, we can see that for any $t\in [0,T]$, 
\begin{align*}
&\bE\left[\int_{0}^{t} \int_{\mathbb{R}^d}\left(\int_{0}^{t} \int_{\mathbb{R}^d} D_{r, z} g(Y_{s-}^n, \eta) \widetilde{N}(\d s, \d \eta)\right)^{2} \nu(\d z) \d r\right]\\
=& \bE\left[\int_{0}^{t} \int_{\mathbb{R}^d}\left(\int_{0}^{t} \int_{\mathbb{R}^d} [g(Y^n_s+D_{r,z}Y^n_{s-}, \eta)-g(Y^n_{s-}, \eta)]\widetilde{N}(\d s, \d \eta)\right)^{2} \nu(\d z) \d r\right]\\
\overset{\eqref{Eq-De-g}}{\leq} &  C \bE \left[\int_{0}^{t} \int_{B_R} \left[ \int_{0}^{t} \int_{B_R}  \|\nabla_y g(\cdot, \eta)\|_\infty^2 |\eta|^{2\mu}  |D_{r,z} Y_{s-}^n|^2  \nu(\d \eta) \d s\right] \nu(\d z) \d r\right]\\
\leq& C \bE \int_0^T\int_{B_R} \left[ \int_0^T  |D_{r,z} Y_{s-}^n|^2 \d s  \right] \nu(\d z) \d r  \leq C \sup_{t\in [0,T]} \|D_{r,z} Y_t^n\|_{L^2(\lambda\times\nu\times\P)}<\infty. 
\end{align*}
Thanks to Theorem 12.15 of \cite{di2009malliavin}, we obtain 
\begin{align*}
&D_{r,z} \int_0^t\!\!\!\int_{B_R} g(Y^n_{s-},\eta)\widetilde{N}(\d s,\d \eta)\\
=&\left[ \int_0^t\!\!\!\int_{B_R} [g(Y^n_{s-}+D_{r,z}Y^n_{s-},\eta)-g(Y^n_{s-},\eta)]\widetilde{N}(\d s,\d\eta)+g(Y^n_{r-},z)\right] \1_{[0,t]}(r). 
\end{align*}
So we obtain $Y^{n+1}_t\in \mathbb{D}_2^1$ and for almost every $r\in [0,t]$,  
\begin{align*}
D_{r,z}Y^{n+1}_t=&g(Y^n_{r-},z)
+\int_r^t [a(Y^n_s+D_{r,z}Y^n_s)-a(Y^n_s)]\d s\\
&+\int_r^t\!\!\!\int_{B_R} [g(Y^n_{s-}+D_{r,z}Y^n_{s-},\eta)-g(Y^n_{s-},\eta)]\widetilde{N}(\d s,\d\eta)
\end{align*}
For any $r\in [0,T]$, denote
$$f^n_r= \E\left[\int_{B_R} \left[ \sup_{r\leq t\leq T} |D_{r,z}Y^n_t|^2 \right]\nu(\d z)\right]$$
Again by Doob's inequality, 
\begin{align*}
f^{n+1}_r =&\E\left[\int_{B_R} \left[ \sup_{r\leq t\leq T}|D_{r,z}Y^{n+1}_t|^2 \right]\nu(\d z)\right]\\
\leq &C\Bigg\{ \int_{B_R} \|g(\cdot,z)\|_\infty^2 \nu(\d z)+ \|\nabla a\|_{\infty}^2 \E\int_r^T\!\!\!\int_{B_R}|D_{r,z}Y^n_{s-}|^2\nu(\d z) \d s\\
&+   \!\int_{B_R}  \E \left[\int_r^T\!\!\!\int_{B_R} \|\nabla_y g(\cdot,\eta)\|_{\infty}^2 | D_{r,z}Y^n_{s-}|^2 \nu(\d \eta) \d s \right] \nu(\d z) \Bigg\}\\
\leq&C\Bigg\{1+T \E\left[\int_{B_R}   \left[ \sup_{r\leq t\leq T}|D_{r,z}Y^n_{t-}|^2 \right]\nu(\d z)\right] \\
&+T\E    \int_{B_R} \left[ \sup_{r\leq t\leq T}|D_{r,z} Y_{t-}^n|^2\right] \nu(\d z) \int_{B_R} |\eta|^{2\mu} \nu(\d \eta)  \Bigg\}\\
=& C+CT f_r^n , 
\end{align*}
where $C$ is independent with $n, r, T$. By choosing $T\leq T_0:= \frac{1}{2C}$, then we have 
$$
f_r^{n}\leq C+\frac{f_r^{n-1}}{2}\leq \cdots\leq 2C + f_r^0=2C + \bE \int_{B_R}| D_{r,z} Y_0 |^2 \nu(\d z). 
$$
Thus, 
\begin{align*}
 \sup_{n\in \N; r\in [0,T] } \E\left[\int_{B_R} \left[ \sup_{ t\in [r,T]} | D_{r,z}Y^n_t|^2\right]\nu(\d z)\right]<\infty. 
\end{align*}
which implies \eqref{DerivateBound} for sufficiently small $T$. For arbitrary $T>T_0$, by the similar argument above, we can see that 
\begin{align*}
&\sup_{n\in \N; r\in [0,T] } \E\left[\int_{B_R} \left[ \sup_{ t\in [r,T]} | D_{r,z}Y^n_t|^2\right]\nu(\d z)\right]\leq  \sup_{n\in \N; r\in [0,T] }f_r^n \\
\leq& 2C+ \sup_{n\in \N; r\in [0,T] }  \bE \int_{B_R}| D_{r,z} Y^n_{T-T_0} |^2 \nu(\d z)\\
\leq& \cdots \leq 2C([T/T_0]+1)+  \sup_{r\in[0,T]}\bE \int_{B_R}| D_{r,z} Y_0 |^2 \nu(\d z)<\infty.  
\end{align*}
So we complete our proof. 

\medspace

(ii)  Choosing $\lambda$ sufficiently large, by \ref{A_nu2} and \eqref{Eq-De-g}, for any $z\in B_{r_0}$, we have 
\begin{align*}
|\nabla_y g(y,z) |\leq&  \|\nabla \si\|_\infty |z|+c(\lambda, \mu) (1-c(\lambda, \mu))^{-1} \big(\|\si\|_\infty^\mu  |z|^\mu+ \|\nabla\si\|_\infty|z| \big)\\
\leq& r_0\|\nabla \si\|_\infty^{-1}+ C\cdot c(\lambda, \mu) <1.  
\end{align*}
This implies that for each $z\in \mathrm{supp}\,\nu\subseteq B_{r_0}$ the map $y\mapsto y+ g(y,z)$ is homeomorphic and $\mI+\nabla_y g(y,z)$ is invertible. Again by \eqref{Eq-De-g}, for any $z\in B_{r_0}$, $\|\nabla_y g(\cdot,z)\|_\infty\leq K(z)\asymp |z|^\mu$. Since $2\mu>\a$, by \eqref{Eq-thetaorder}, 
\begin{align*}
\int_{B_{r_0} } K(z)^2\nu(\d z)\leq & C \int_{B_{r_0}} |z|^{2\mu} \nu_2(\d z)\leq C r_0^{2\mu-\a}<\infty. 
\end{align*}
Noting that $\sigma\in C_b^{1+\delta}$, by \eqref{Eq-de-g} and the regularity estimates for $u$, one can also check that 
$$
|\nabla_y g(y,z)-\nabla_y g(y',z)|\leq L(z) |y-y'|^{\delta\wedge \mu}
$$ 
and  $L(z)\asymp |z|^\mu$. So we also have 
$$
\int_{B_{r_0}} L(z)^2\nu(\d z)<\infty.
$$
Thanks to \cite[Theorem 3.11]{kunita2004stochastic}, $\{Y_t(x)\}$ defines a $C^1$-stochastic flow. So  $\{X_t(x)\}$ also forms a $C^1$-stochastic flow.  
\end{proof}

 In \cite{davie2007uniqueness}, Davie  proved a remarkable result, which says that given a bounded Borel measurable mapping $b$, for almost all Brownian paths $W_t(\om)$, there is a unique $\theta_t(\om)$ satisfying the random ODE: $\d \theta_t(\om) =b(\theta_t(\om)+W_t(\om))\d t$. This type of uniqueness for SDEs with jumps was first studied by Priola in \cite{priola2018davie}. In the light of \cite[Theorem 1.1]{priola2018davie}, we show a path-by-path uniqueness result for \eqref{Eq-SDE} below. 
\bp\label{Cor-Main}
Suppose $Z$ is a non-degenerate symmetric $\a$-stable process, $\sigma=\mathbb{I}$ and  
$b\in C_b^\beta$ with $\beta\in (1-\frac{\a}{2},1)$. Then there is a $\bP$-null set $\cN\subseteq \Om$ such that for any $\omega\notin \cN$, the following ODE: 
\be\label{Eq-RODE}
\frac{\d z_t(\om)}{\d t}=b\left(z_t(\omega)+Z_t(\om)\right), \quad z_0=x 
\ee
admits a unique solution in $C(\mR_+; \mR^d)$. 
\ep
\begin{proof}
 Fix $\gamma\in (1-\frac{\a}{2}, \beta)$. By Remark \ref{Rek-Const} (2) and interpolation theorem,  there is a constant $\lambda>$ such that the following resolvent equation  
$$
\lambda u_\lambda- L u_\lambda-b\cdot \nabla u_\lambda=b 
$$
has a solution $u_\lambda\in C_b^{\a+\gamma}$ satisfying $\|\nabla u_\lambda\|_\infty<1/3$. Thanks to \cite[Theorem 6.6 and Theorem 1.1]{priola2018davie}, \eqref{Eq-RODE} has exactly one solution $z$ in $C(\mR_+; \mR^d)$. 
\end{proof}

\section*{Acknowledgements}
The author is very grateful to Professor Enrico Priola who read the first draft and pointed out some mistakes and also to Professor Moritz Kassmann and Professor Xicheng Zhang for their valuable discussion.

\bibliographystyle{alpha}

\bibliography{mybib}

\newcommand{\etalchar}[1]{$^{#1}$}
\begin{thebibliography}{MPMBN{\etalchar{+}}13}

\bibitem[BCD11]{bahouri2011fourier}
Hajer Bahouri, Jean-Yves Chemin, and Rapha{\"e}l Danchin.
\newblock {\em Fourier analysis and nonlinear partial differential equations},
  volume 343.
\newblock Springer Science \& Business Media, 2011.

\bibitem[BK05]{bass2005holder}
Richard~F Bass and Moritz Kassmann.
\newblock H{\"o}lder continuity of harmonic functions with respect to operators
  of variable order.
\newblock {\em Communications in Partial Differential Equations},
  30(8):1249--1259, 2005.

\bibitem[BL02a]{bass2002transition}
Richard Bass and David Levin.
\newblock Transition probabilities for symmetric jump processes.
\newblock {\em Transactions of the American Mathematical Society},
  354(7):2933--2953, 2002.

\bibitem[BL02b]{bass2002harnack}
Richard~F Bass and David~A Levin.
\newblock Harnack inequalities for jump processes.
\newblock {\em Potential Analysis}, 17(4):375--388, 2002.

\bibitem[BSW13]{bottcher2013levy}
Bj{\"o}rn B{\"o}ttcher, Ren{\'e} Schilling, and Jian Wang.
\newblock {\em L{\'e}vy matters. III. L\'evy-Type Processes: Construction,
  Approximation and Sample Path Properties}, volume 2099.
\newblock Springer, 2013.

\bibitem[CCV11]{caffarelli2011regularity}
Luis Caffarelli, Chi~Hin Chan, and Alexis Vasseur.
\newblock Regularity theory for parabolic nonlinear integral operators.
\newblock {\em Journal of the American Mathematical Society}, 24(3):849--869,
  2011.

\bibitem[CK03]{chen2003heat}
Zhen-Qing Chen and Takashi Kumagai.
\newblock Heat kernel estimates for stable-like processes on d-sets.
\newblock {\em Stochastic Processes and their applications}, 108(1):27--62,
  2003.

\bibitem[CMZ07]{chen2007new}
Qionglei Chen, Changxing Miao, and Zhifei Zhang.
\newblock A new bernstein’s inequality and the 2d dissipative
  quasi-geostrophic equation.
\newblock {\em Communications in mathematical physics}, 271(3):821--838, 2007.

\bibitem[CSZ18]{chen2018stochastic}
{Z}hen-{Q}ing {C}hen, {R}enming {S}ong, and {X}icheng {Z}hang.
\newblock {S}tochastic flows for {L}\'evy processes with {H}\"{o}lder drifts.
\newblock {\em Revista Matemática Iberoamericana}, 34(4), 2018.

\bibitem[CZZ21]{chen2021supercritical}
Zhen-Qing Chen, Xicheng Zhang, and Guohuan Zhao.
\newblock Supercritical {{SDE}}s driven by multiplicative stable-like
  {L}{\'e}vy processes.
\newblock {\em To apear in {T}ransaction of {A}merican {M}athematical
  {S}ociety.}, 2021.

\bibitem[Dav07]{davie2007uniqueness}
Alexander~M Davie.
\newblock Uniqueness of solutions of stochastic differential equations.
\newblock {\em International Mathematics Research Notices}, 2007, 2007.

\bibitem[DJZ18]{dong2018dini}
Hongjie Dong, Tianling Jin, and Hong Zhang.
\newblock Dini and schauder estimates for nonlocal fully nonlinear parabolic
  equations with drifts.
\newblock {\em Analysis \& PDE}, 11(6):1487--1534, 2018.

\bibitem[DK13]{dong2013schauder}
Hongjie Dong and Doyoon Kim.
\newblock Schauder estimates for a class of non-local elliptic equations.
\newblock {\em Discrete Contin. Dyn. Syst.}, 33(6):2319–2347, 2013.

\bibitem[DN{\O}P09]{di2009malliavin}
Giulia Di~Nunno, Bernt~Karsten {\O}ksendal, and Frank Proske.
\newblock {\em {M}alliavin calculus for {L}\'evy processes with applications to
  finance}, volume~2.
\newblock Springer, 2009.

\bibitem[dRMP20]{de2020schauder}
Paul-{\'E}ric~Chaudru de~Raynal, Stephane Menozzi, and Enrico Priola.
\newblock Schauder estimates for drifted fractional operators in the
  supercritical case.
\newblock {\em Journal of Functional Analysis}, 278(8):108425, 2020.

\bibitem[HL18]{huang2018euler}
Xing Huang and Zhong-Wei Liao.
\newblock The {E}uler-{M}aruyama method for {{SDE}}s with {H}{\"o}lder drift
  and $\alpha$-stable noise.
\newblock {\em Stochastic Analysis and Applications}, 36(1):28--39, 2018.

\bibitem[HP14]{haadem2014construction}
Sven Haadem and Frank Proske.
\newblock On the construction and {M}alliavin differentiability of solutions of
  {L}\'evy noise driven {{SDE}}s with singular coefficients.
\newblock {\em Journal of Functional Analysis}, 266(8):5321--5359, 2014.

\bibitem[HWW20]{hao2020schauder}
Zimo Hao, Zhen Wang, and Mingyan Wu.
\newblock Schauder's estimates for nonlocal equations with singular {L}\'evy
  measures.
\newblock {\em arXiv preprint arXiv:2002.09887}, 2020.

\bibitem[Jac05]{jacob2005pseudo}
Niels Jacob.
\newblock {\em Pseudo differential operators and Markov processes}, volume
  I-III.
\newblock Imperial College Press, 2001--2005.

\bibitem[JS01]{jacob2001levy}
Niels Jacob and Ren{\'e}~L Schilling.
\newblock L{\'e}vy-type processes and pseudodifferential operators.
\newblock In {\em L{\'e}vy processes: Theory and Applications}, pages 139--168.
  Springer, 2001.

\bibitem[Kas09]{kassmann2009priori}
Moritz Kassmann.
\newblock A priori estimates for integro-differential operators with measurable
  kernels.
\newblock {\em Calculus of Variations and Partial Differential Equations},
  34(1):1--21, 2009.

\bibitem[KR05]{krylov2005strong}
{N}icolai~{V} {K}rylov and {M}ichael {R}\"ockner.
\newblock {S}trong solutions of stochastic equations with singular time
  dependent drift.
\newblock {\em {P}robability Theory and Related Fields}, 131(2):154--196, 2005.

\bibitem[KS19]{kuhn2019strong}
Franziska K{\"u}hn and Ren{\'e}~L Schilling.
\newblock Strong convergence of the {E}uler-{M}aruyama approximation for a
  class of {L}{\'e}vy-driven {{SDE}}s.
\newblock {\em Stochastic Processes and their Applications}, 129(8):2654--2680,
  2019.

\bibitem[Kul19]{kulik2019weak}
Alexei~M Kulik.
\newblock On weak uniqueness and distributional properties of a solution to an
  {SDE} with $\alpha$-stable noise.
\newblock {\em Stochastic Processes and their Applications}, 129(2):473--506,
  2019.

\bibitem[Kun04]{kunita2004stochastic}
Hiroshi Kunita.
\newblock Stochastic differential equations based on {L}\'evy processes and
  stochastic flows of diffeomorphisms.
\newblock In {\em Real and stochastic analysis}, pages 305--373. Springer,
  2004.

\bibitem[LZ19]{lin2019nonlocal}
Chengcheng Lin and Guohuan Zhao.
\newblock Nonlocal elliptic equation in {H}\" older space and the martingale
  problem.
\newblock {\em arXiv preprint arXiv:1907.00588}, 2019.

\bibitem[MBP10]{meyer2010construction}
Thilo Meyer-Brandis and Frank Proske.
\newblock Construction of strong solutions of {{SDE}}'s via {M}alliavin
  calculus.
\newblock {\em Journal of Functional Analysis}, 258(11):3922--3953, 2010.

\bibitem[MNP15]{mohammed2015sobolev}
Salah-Eldin~A Mohammed, Torstein~K Nilssen, and Frank~N Proske.
\newblock {S}obolev differentiable stochastic flows for {{SDE}}s with singular
  coefficients: Applications to the transport equation.
\newblock {\em The Annals of Probability}, 43(3):1535--1576, 2015.

\bibitem[MP92]{mikulevivcius1992cauchy}
Remigijus Mikulevicius and Henrikas Pragarauskas.
\newblock On the {C}auchy problem for certain integro-differential operators in
  {S}obolev and {H}{\"o}lder spaces.
\newblock {\em Lithuanian Mathematical Journal}, 32(2):238--264, 1992.

\bibitem[MPMBN{\etalchar{+}}13]{menoukeu2013variational}
Olivier Menoukeu-Pamen, Thilo Meyer-Brandis, Torstein Nilssen, Frank Proske,
  and Tusheng Zhang.
\newblock A variational approach to the construction and {M}alliavin
  differentiability of strong solutions of {{SDE}}’s.
\newblock {\em Mathematische Annalen}, 357(2):761--799, 2013.

\bibitem[MX18]{mikulevivcius2018rate}
R~Mikulevi{{c}}ius and Fanhui Xu.
\newblock On the rate of convergence of strong {E}uler approximation for
  {{SDE}}s driven by {L}{\'e}vy processes.
\newblock {\em Stochastics}, 90(4):569--604, 2018.

\bibitem[{P}ri12]{priola2012pathwise}
{E}nrico {P}riola.
\newblock {P}athwise uniqueness for singular {{S}{D}{E}}s driven by stable
  processes.
\newblock {\em {O}saka {J}ournal of {M}athematics}, 49(2):421--447, 2012.

\bibitem[{P}ri15]{priola2015stochastic}
{E}nrico {P}riola.
\newblock Stochastic flow for {{SDE}}s with jumps and irregular drift term.
\newblock {\em Banach Center Publications}, 105(1):421--44193--210, 2015.

\bibitem[Pri18]{priola2018davie}
Enrico Priola.
\newblock Davie's type uniqueness for a class of {{SDE}}s with jumps.
\newblock {\em Annales de l'Institut Henri Poincar{\'e}, Probabilit{\'e}s et
  Statistiques}, 54(2):694--725, 2018.

\bibitem[Sil06]{silvestre2006holder}
Luis Silvestre.
\newblock H{\"o}lder estimates for solutions of integro-differential equations
  like the fractional laplace.
\newblock {\em Indiana University mathematics journal}, pages 1155--1174, 2006.

\bibitem[{S}il12]{silvestre2012differentiability}
{L}uis {S}ilvestre.
\newblock {O}n the differentiability of the solution to an equation with drift
  and fractional diffusion.
\newblock {\em Indiana University Mathematics Journal}, 61(2):557--584, 2012.

\bibitem[SS10]{schilling2009symbol}
Ren{\'e}~L Schilling and Alexander Schnurr.
\newblock The symbol associated with the solution of a stochastic differential
  equation.
\newblock {\em Electronic Journal of Probability}, 15:1369--1393, 2010.

\bibitem[Tri92]{triebel1992theory}
Hans. Triebel.
\newblock {\em Theory of function spaces II}, volume~84.
\newblock Birkh\"auser Basel, 1992.

\bibitem[TTW74]{tanaka1974perturbation}
{H}iroshi {T}anaka, {M}asaaki {T}suchiya, and {S}hinzo {W}atanabe.
\newblock {P}erturbation of drift-type for {L}{\'e}vy processes.
\newblock {\em {J}ournal of {M}athematics of {K}yoto {U}niversity},
  14(1):73--92, 1974.

\bibitem[Ver80]{veretennikov1980strong2}
Alexander~Yur'evich Veretennikov.
\newblock On strong solutions and explicit formulas for solutions of stochastic
  integral equations.
\newblock {\em Matematicheskii Sbornik}, 153(3):434--452, 1980.

\bibitem[XXZZ20]{xia2020lqlp}
Pengcheng Xia, Longjie Xie, Xicheng Zhang, and Guohuan Zhao.
\newblock ${L}^q({L}^p)$-theory of stochastic differential equations.
\newblock {\em Stochastic Processes and their Applications}, 130(8):5188--5211,
  2020.

\bibitem[{Z}ha11]{zhang2011stochastic}
{X}icheng {Z}hang.
\newblock {S}tochastic homeomorphism flows of {S}{D}{E}s with singular drifts
  and {S}obolev diffusion coefficients.
\newblock {\em {E}lectronic {J}ournal of {P}robability}, 16:1096--1116, 2011.

\bibitem[Zha13]{zhang2013stochastic}
Xicheng Zhang.
\newblock Stochastic differential equations with {S}obolev drifts and driven by
  $\alpha$-stable processes.
\newblock {\em Annales de l'IHP Probabilit{\'e}s et statistiques},
  49(4):1057--1079, 2013.

\bibitem[{Z}ha16]{zhang2016stochastic}
{X}icheng {Z}hang.
\newblock {S}tochastic differential equations with {S}obolev diffusion and
  singular drift and applications.
\newblock {\em {T}he {A}nnals of {A}pplied {P}robability}, 26(5):2697--2732,
  2016.

\bibitem[Zha19]{zhao2019weak}
Guohuan Zhao.
\newblock Weak uniqueness for {{SDE}}s driven by supercritical stable processes
  with {H}{\"o}lder drifts.
\newblock {\em Proceedings of the American Mathematical Society},
  147(2):849--860, 2019.

\bibitem[Zvo74]{zvonkin1974transformation}
Alexander~K Zvonkin.
\newblock A transformation of the phase space of a diffusion process that
  removes the drift.
\newblock {\em Mathematics of the USSR-Sbornik}, 22(1):129, 1974.

\bibitem[ZZ18]{zhang2018dirichlet}
Xicheng Zhang and Guohuan Zhao.
\newblock Dirichlet problem for supercritical nonlocal operators.
\newblock {\em arXiv preprint arXiv:1809.05712}, 2018.

\bibitem[ZZ21]{zhang2020stochastic}
Xicheng Zhang and Guohuan Zhao.
\newblock Stochastic lagrangian path for leray’s solutions of 3d
  navier--stokes equations.
\newblock {\em Communications in Mathematical Physics}, 381(2):491--525, 2021.

\end{thebibliography}

\end{document}